\documentclass[arma]{svjour}
\usepackage{graphics,color}
\usepackage{amssymb}
\usepackage{latexsym}
\usepackage{here}

\makeatletter\def\extrarowheight#1{\noalign{\@tempdima\ht\@arstrutbox
\advance\@tempdima#1\ht\@arstrutbox\@tempdima}}\makeatother

\DeclareMathAlphabet{\mathbi}{OML}{cmm}{b}{it} 
\newcommand{\bel}{\begin{equation}\label}
\newcommand{\ee}{\end{equation}}
\newcommand{\beq}{\begin{eqnarray}\label}
\newcommand{\eq}{\end{eqnarray}}
\newcommand{\bc}{\begin{center}}
\newcommand{\ec}{\end{center}}
\newcommand{\bu}{\mathbi{u}}
\newcommand{\bdf}{\mathbi{f}} 
\newcommand{\bx}{\mathbi{x}}
\newcommand{\bk}{\mathbi{k}}
\newcommand{\cE}{\mathcal{E}}
\newcommand{\bit}{\begin{itemize}}
\newcommand{\eit}{\end{itemize}}
\newcommand{\ben}{\begin{enumerate}}
\newcommand{\een}{\end{enumerate}} 
\newcommand{\bw}{\mbox{\boldmath$\omega$}}

\newcommand{\kn}{\kappa_{n}}
\newcommand{\knp}{\kappa_{n+1}}

\newcommand{\knr}{\kappa_{n,r}}

\newcommand{\tplus}{(\Delta t^{i,j}_{+})}
\newcommand{\tij}{t_{0}^{i,j}}
\newcommand{\tg}{(\Delta t)_{g}}
\newcommand{\tb}{(\Delta t)_{b}}

\newcommand{\Idt}{\int_{\Delta t}}

\newcommand{\calE}{\mathcal{E}}
\newcommand\shalf{\ensuremath{{\scriptstyle\frac{1}{2}}}}
\newcommand\sixth{\ensuremath{{\scriptstyle\frac{1}{6}}}}

\newcommand\squart{\ensuremath{{\textstyle\frac{1}{4}}}}

\newcommand\thfth{\ensuremath{{\textstyle\frac{3}{5}}}}

\newcommand{\enk}{e^{-\nu k_{1}^{2}\Delta t}}

\newcommand{\nk}{\nu k_{1}^{2}}

\newcommand{\Rey}{\mathrm{Re}}
\newcommand{\Gr}{\mathrm{Gr}}
\newcommand{\mpo}{(1+\mu)}
\newcommand{\lbd}{\lambda}
\newcommand{\X}{\rule[0.0cm]{0.35cm}{0.1cm}} 
\newcommand{\Q}{\rule[0.0cm]{0.35cm}{0.0cm}} 
\newcommand{\Y}{\rule[0.0cm]{0.35cm}{0.25cm}} 

\begin{document}
\title{Intermittency and regularity issues in $3D$ Navier-Stokes 
turbulence\footnote{Original version Nov 03; 1st revised-version 
6th June 04}}
\author{J. D. Gibbon and Charles R. Doering
}                     
\maketitle
\begin{abstract}
Two related open problems in the theory of $3D$ Navier-Stokes turbulence are
discussed in this paper. The first is the phenomenon of intermittency in the
dissipation field. Dissipation-range intermittency was first discovered
experimentally by Batchelor and Townsend over fifty years ago. It is
characterized by spatio-temporal binary behaviour in which long, quiescent
periods in the velocity signal are interrupted by short, active `events'
during which there are violent fluctuations away from the average. The second
and related problem is whether solutions of the $3D$ Navier-Stokes equations
develop finite time singularities during these events. This paper shows that
Leray's weak solutions of the three-dimensional incompressible Navier-Stokes
equations can have a binary character in time. The time-axis is split into
`good' and `bad' intervals: on the `good' intervals solutions are bounded and
regular, whereas singularities are still possible within the `bad' intervals.
An estimate for the width of the latter is very small and decreases with
increasing Reynolds number. It also decreases relative to the lengths of the
good intervals as the Reynolds number increases. Within these `bad' intervals,
lower bounds on the local energy dissipation rate and other quantities, such
as $\|\bu(\cdot,\,t)\|_{\infty}$ and $\|\nabla\bu(\cdot,\,t)\|_{\infty}$, are
very large, resulting in strong dynamics at sub-Kolmogorov scales.
Intersections of bad intervals for $n\geq 1$ are related to Scheffer's
potentially singular set in time. It is also proved that the Navier-Stokes
equations are conditionally regular provided, in a given `bad' interval, the
energy has a lower bound that is decaying exponentially in time.
\end{abstract}


\section{Introduction}\label{intro}

\begin{figure}
\resizebox{0.75\textwidth}{!}{%
\includegraphics{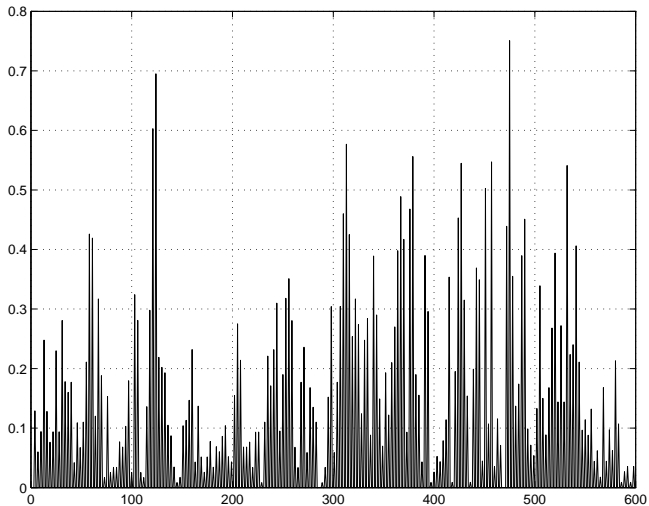}
} 
\caption{A typical example of dissipation-range intermittency from wind tunnel 
turbulence where hot wire anemometry has been used to measure the longitudinal 
velocity derivative at a single point (D. Hurst and J. C. Vassilicos). The 
horizontal axis spans 8 integral time scales. The Taylor micro-scale based 
Reynolds number is about 200.}
\label{jcvfig1}
\end{figure}

The questions to be addressed in this paper concern the nature and behaviour
of intermittent high Reynolds number solutions of the three-dimensional
Navier-Stokes equations. Relatively quiescent flows can exist in nature at
high Reynolds numbers. There are also flows that appear turbulent for all
practical purposes but are nevertheless smooth at appropriately small length
scales. Typically turbulent high Reynolds number Navier-Stokes flows, however,
generally display a specific hallmark which is called dissipation-range
intermittency. This was first discovered by Batchelor and Townsend \cite{BT49}
and manifests itself in violent fluctuations of very short duration in the
energy dissipation rate. These fluctuations away from the average are
interspersed by quieter, longer periods in the dynamics. The data in Figure
\ref{jcvfig1} is an illustration of a typically intermittent signal
representing a velocity derivative versus time recorded at a single point in
space.
For theoretical studies of the Navier-Stokes equations it is usual to express 
energy dissipation in the $L^{2}$-volume-integrated sense. Postponing the full 
definition of system variables, if the global energy $H_{0}$ and 
enstrophy $H_{1}$ for a three-dimensional Navier-Stokes velocity field 
$\bu(\bx,\,t)$ are defined as
\bel{int1}
H_{0}(t) = \int_{V}|\bu|^2\,dV\hspace{2cm}H_{1}(t) = \int_{V}|\nabla\bu|^2\,dV
\ee
then it is well known that while the energy $H_{0}$ is a uniformly bounded 
function of time, only the long-time averaged energy dissipation rate 
$\epsilon_{av} = \nu L^{-3} \left<H_{1}\right>$ is known to be bounded 
\cite{Leray}, whereas the behaviour of $\epsilon (t) = \nu L^{-3}H_{1}(t)$ 
pointwise in time may be wildly fluctuating or even singular. The 
long-time average $\left<\cdot\right>$ is defined later in this section.

Whether $H_{1}(t)$ becomes singular in a finite time is an open question
intimately related to the Navier-Stokes regularity problem. So long as
$H_{1}(t)$ is finite the solution is smooth and unique; any finite time
singularity must be accompanied by a divergence of $H_{1}(t)$ at that time.
From Leray onwards \cite{Leray}, this question has led to a long and rich
literature on the nature of weak solutions
\cite{Lady,Serrin,Foiasstat1,Foiasstat3,FGT,Temam,CF,FMRT,MB}. In terms of
physical length scales, the boundedness of $\left<H_{1}\right>$ allows the
time averaged turbulent energy dissipation rate per unit volume $\epsilon_{av}
= \nu L^{-3}\left<H_{1}\right>$ to be used naturally in forming the inverse
Kolmogorov length
\bel{Koldef1}
\eta_{K}^{-1} = \left(\epsilon_{av}/ \nu^{3}\right)^{1/4}.
\ee
Spiky behaviour in $H_{1}(t)$, causing loss of resolution in large-scale
computations, could mean that significant energy lies in wave-numbers 
$k > \eta_{K}^{-1}$ in the dissipation range of the energy spectrum. 
The ubiquity of dissipation-range intermittency in turbulent flows suggests 
that it should occur naturally in mathematical 
analyses of the Navier-Stokes equations. Strong temporal excursions in
$H_{1}(t)$ are clear candidates for the formation of singularities and may be
related to potentially singular solutions of the three-dimensional Euler
equations \cite{Kerr,BKM,CFM,MB}, although nothing has been rigorously proved
in this respect. Beginning with Leray's seminal paper in 1934, rigorous
methods of analysis on the full three-dimensional domain have led to seventy
years of literature on the Navier-Stokes regularity problem but have yet to
settle this question definitively. Short-time regularity has also been known
for many years, as have various interesting partial and conditional regularity
results
\cite{Lady,Serrin,Foiasstat1,Foiasstat3,FGT,Temam,CF,FMRT,MB,Scheffer76a,CKN}.

Batchelor and Townsend \cite{BT49} also suggested that energy dissipation is 
not distributed evenly across the
full three-dimensional spatial domain but is clustered into smaller spots in
the flow with the energy associated with the small-scale components being
distributed unevenly in space and roughly confined to regions which become
smaller with eddy size \cite{Kuo71}. In contrast to Kolmogorov's traditional
theory that implies that energy dissipation is space-filling \cite{Frisch},
Mandelbrot suggested that the spatial set on which it occurs is actually
fractal \cite{Man1}. In experimental investigations of the energy dissipation
rate in several laboratory flows, and in the atmospheric surface layer,
Sreenivasan \& Meneveau \cite{SM88,MS91} interpreted the evident intermittent
nature of their signals in terms of multifractals (see also \cite{Sreenirev}).
Zeff \textit{et al} \cite{Lathrop} have shown how more recent technical
advances have made it possible to measure each derivative of all three
velocity components to obtain a fuller experimental picture of the energy
dissipation at a point in a flow. In general, numerical simulations and
experiments suggest that respectively quasi-one and two-dimensional tubes and
sheets are the favoured low-dimensional sets on which vorticity and strain
appear to accumulate \cite{HY,DCB,VM94} although, as Galanti \& Tsinober
\cite{GT} and Tsinober \cite{Tsinjcv} have pointed out, there are significant
differences between these two sets. It is also true that explaining these
phenomena in the simple geometrical terms of tubes and sheets is a visual
over-simplification of much more complicated dynamical spatial structures at
small scales. Examples of this are the spiral vortex structures introduced 
by Lundgren \cite{Ld1} and discussed in detail by Vassilicos \& Hunt \cite{VH},
Flohr \& Vassilicos \cite{FV} and Angilella \& Vassilicos \cite{AV1,AV2}.

Dissipation-range intermittency is a well established, experimentally
observable phenomenon; its appearance in systems other than the Navier-Stokes
equations has been discussed in an early and easily accessible paper by Frisch
\& Morf \cite{FM}. One symptom of its occurrence is the deviation of the
`flatness' of a velocity signal (the ratio of the 4th order moment to the
square of the second order moment) from the value of 3 that usually holds for
Gaussian statistics \cite{Homoturb}. More subtle is the phenomenon of
inertial-range intermittency that has exercised the ingenuity of that section
of the physics community that focuses on scaling methods. In the inertial
range, Kolmogorov's theory predicts that the exponent, $\zeta_{p}$, of the
$p$th velocity structure function should vary linearly with $p$, whereas
experimental data shows that the $(\zeta_{p},\,p)$ relation is a concave curve
lying below the line $p/3$ for $p\geq 3$. This departure from Kolmogorov
scaling in the inertial range, and therefore from the five-thirds law, is
termed inertial-range intermittency. An extensive literature is quoted in
Frisch's book \cite{Frisch}. Studies in weak turbulence, applicable to
predominantly dispersive systems, have been pioneered by Zakharov. These ideas
can be found in papers by Zakharov, L'vov and Falkovich \cite{Z1} and Zakharov
\cite{Z2}.

To prove that solutions of the Navier-Stokes equations (usually taken in a
periodic box \cite{Temam,CF,FMRT,DGbook}) are typically intermittent in
space-time poses formidable technical challenges to the mathematician.
Analysis on time-evolving fractal domains with ill-defined boundary conditions
is not advanced enough to gain rigorous results by concentrating on one
fractal `spot' in the flow. The partial regularity result of Scheffer
\cite{Scheffer76a} proving that the potentially singular set in time alone has
zero half-dimensional Hausdorff measure has been of great influence (see also
\cite{FT} for a short proof). Following this, the most significant space-time
result has been that of Caffarelli, Kohn \& Nirenberg \cite{CKN} who showed
that the potentially singular set in space-time has zero one-dimensional
Hausdorff measure. This implies that if singularities do exist they must be
relatively rare. Lin \cite{Lin} and Choe \& Lewis \cite{CL} have recently
provided shorter proofs of this result.

The more realistic option adopted here is to show that solutions of the
Navier-Stokes equations can have a \textit{binary nature in time} in which the
time-axis is divided into what are designated as \textit{good} and
\textit{bad} intervals. On the good intervals the Navier-Stokes equations are
uniformly regular. The bad intervals are shown to be very small in width with
an upper bound that decreases with increasing Reynolds number and which also
decreases relative to the widths of the good intervals. Within the bad
intervals very large lower bounds are shown to exist on both the local-in-time
energy dissipation rate and several other quantities, such as
$\|\bu(\cdot,\,t)\|_{\infty}$ and $\|\nabla\bu(\cdot,\,t)\|_{\infty}$. The
corresponding local length scales within these intervals are, \textit{at
best}, comparatively much smaller than the Kolmogorov length. The regularity
question within the bad intervals is still open, so only weak solutions are
known to exist there. The great difficulties encountered by computational
fluid dynamicists in resolving turbulent flows even for modestly high Reynolds
numbers could be because of this binary behaviour.

These results, which are summarized in \S2, have been obtained through the use
of a set of quantities $\kn(t)$ that have been introduced in previous papers
\cite{DGbook,BDGM,Doering02}. The more physical aspects of these ideas have
been laid out previously in a short paper \cite{GD03}; the present paper gives
a detailed and more advanced account of the methods and results reported
there. The $\kn(t)$ have the dimensions of inverse lengths and are formed from
ratios of $L^2$-norms of derivatives of the velocity field. Together with the
periodicity of the domain, the $L^{2}$-spatial integration within the $\kn(t)$
means that spatially intermittent effects are included implicitly and cannot
be averaged away. Clearly they are not the same quantities as those measured
by experimentalists, such as the energy dissipation rate, but estimates for
them are rigorous, making no appeal to any approximations, and ultimately lead
to information on the energy dissipation. Physically they can be considered as
a measure of the $2n$th moment of the energy spectrum. Their time-dependence
is explicit and their long-time averages $\left<\kn\right>$ are uniformly
bounded \cite{Doering02}. Pointwise in time their binary nature appears for
each value of $n\geq 2$.

The phrase `solutions can have a binary nature in time' has been used above in
the following sense: if no bad intervals occur, then solutions of the
Navier-Stokes equations are bounded for all time. In this sense, the results
in this paper are different from conventional short-time regularity proofs
because loss of regularity can only occur in the bad intervals. This is
consistent with Scheffer's partial regularity result \cite{Scheffer76a}: if
singularities exist at points in time then these must be clustered within the
intersection of the bad intervals. They are also consistent with the well
known problems of computational resolution in three-dimensional turbulence:
the bounds indicate a structure so fine that it would be extremely difficult
to distinguish between regular and singular solutions, despite the possibility
of a lack of sharpness.
\begin{table}[ht]
\bc
\caption{\label{tab1}Definitions of the main parameters in the paper.}
\begin{tabular}{lll} 
\hline\noalign{\smallskip}
\extrarowheight{9pt}
Quantity & Definition & Comment\\
\noalign{\smallskip}\hline\noalign{\smallskip}
Box length & $L$ & \\
Forcing length scale & $\ell$ & $\ell \leq L/2\pi$\\
Average forcing & $f^{2} = L^{-3}\|\bdf\|_{2}^{2}$ & \\
Narrow-band forcing & $\|\bdf\|_{2}^{2} \approx \ell^{2n}\|\nabla^{n}\bdf\|_{2}^{2}$&\\
Average velocity & $U^{2}=L^{-3}\left<\|\bu\|_{2}^{2}\right>$ & \\
Grashof No & $\Gr = f\ell^{3}\nu^{-2}$ & \\
Reynolds No \cite{DF} & $\Rey = U\ell\nu^{-1}$ &  $\Gr^{1/2} 
\leq c\,\Rey~\mbox{as}~\Gr\to\infty$\\
\noalign{\smallskip}\hline
\end{tabular}
\ec
\end{table}
At this point it is appropriate to describe the Navier-Stokes system of
partial differential equations considered on a periodic cube $V = [0,\,L]^{3}$
\bel{I1}
\bu_{t}+\bu\cdot\nabla\bu = \nu\Delta\bu - \nabla p + \bdf(\bx),
\hspace{1.5cm}\nabla\cdot\bu = 0
\ee
where $\nu$ is the kinematic viscosity and $p$ is the pressure. The applied
body force $\bdf(\bx)$ is taken to be mean-zero and divergence-free so,
without loss of generality, the solution $\bu(\bx,t)$ is mean-zero at all
times. For simplicity narrow-band body forces with a single length scale
$\ell$ are considered; that is, with Fourier components only at wave-number $k
= \ell^{-1}$ and $\ell \leq L/2\pi$. For finite energy initial data the
Navier-Stokes equations admit weak solutions in $L^{2}(V)$ at each instant of
time, with finite time integrals of the $L^{2}$-norms of the velocity
gradients. With the assumption of narrow-band forcing, norms of gradients of
$\bdf(\bx)$ are related to the norm of $\bdf(\bx)$; these can all be found in
Table \ref{tab1}. Also found there are the definitions of the root-mean-square
velocity scale $U$ and the Reynolds number $\Rey$. The angled brackets
$\left<\cdot\right>$ denote the long-time average
\bel{I8}
\left<\Phi(\cdot) \right> = \mbox{\textsc{lim}}_{t\to\infty}\left(\frac{1}{t}
\int_{0}^{t}\Phi(s)\,ds\right).
\ee
$\mbox{\textsc{lim}}$ is a generalized long-time limit for functionals of
(weak) statistical solutions of the Navier-Stokes equations 
\cite{Foiasstat1,Foiasstat3}. The square of the $L^2$-norm $\|\cdot\|_{2}^{2}$ 
is defined as
\bel{L2def}
\|\bdf\|_{2}^{2} = \int_{V}|\bdf|^{2}\,dV
\ee
The Grashof number $\Gr$ is the natural control parameter, not the Reynolds 
number $\Rey$, but it is clear that high Reynolds number solutions may be 
achieved if the Grashof number $\Gr$ is sufficiently high; indeed, Doering 
\& Foias \cite{DF} have proved that for body-forced Navier-Stokes flows 
such as these $\Gr^{1/2}\leq c\,\Rey$ as $\Gr\to \infty$. 
Next are defined the time dependent quantities 
\bel{Hndef}
H_{n}(t) = \int_{V}|\nabla^{n}\bu(\bx,t)|^{2}\,dV = 
\sum_{\bk}k^{2n}|\hat{\bu}(\bk,t)|^{2}.
\ee
where $\hat{\bu}(\bk,t)$ is the Fourier transform of $\bu(\bx,t)$. 
For higher derivatives of $\bu$ and $\bdf$, the important quantities that will
be used in this paper are
\bel{I3}
F_{n} = \int_{V} \left(|\nabla^{n}\bu|^{2} + 
\tau^{2}|\nabla^{n}\bdf|^{2}\right)\,dV
\ee
where $\tau$ is a characteristic time the origin of which is discussed later 
in \S\ref{standard}. In approaching this problem conventional $L^2$-norms 
have been used in (\ref{I3}) to avoid difficulties with the pressure field. 

Not all interesting small-scale spatial behaviour can be averaged away on a 
finite periodic box; spatial events must show up in some temporal manner. The 
next step is to define the quantities
\bel{I4}
\knr(t) = \left(F_{n}/F_{r}\right)^{1/2(n-r)}
\ee
The particular one of interest is $\kn = \kappa_{n,0}$ where the $r=0$ label
has been dropped for convenience in the first. $\kn^{2n}$ can be interpreted as
being related to the $2n$th moment of the energy spectrum \cite{Doering02}; 
this can be seen from writing $|v|^{2} = |\hat{\bu}|^{2} + |\hat{\bdf}|^{2}$ 
and then applying Parseval's equality to (\ref{I4}) (with $r=0$)
\bel{I5a}
\left[\kn(t)\right]^{2n} 
= \frac{\sum_{k}k^{2n}|v(\bk,t)|^{2}}{\sum_{k}|v(\bk,t)|^{2}}.
\ee
Moreover, the $\kn$ are ordered in magnitude for all $t\geq 0$
\bel{I5}
L^{-1} \leq \kappa_{1} \ldots \leq \kn \leq \knp \leq \dots
\ee
which is simply a result of H\"{o}lder's inequality. The full $\kappa_{n,r}$ 
are also ordered such that $\knr \leq \kappa_{n,r+1}$ for $r+1 < n$.
\begin{table}[ht]
\bc
\caption{\label{tab2}Definitions of the main variables and constants 
in the paper. The parameter $\delta$, which lies in the range 
$0<\delta<\sixth$, is ignored hereafter.}
\begin{tabular}{lll}
\hline\noalign{\smallskip}
\extrarowheight{9pt}
 & Definition & \\
\noalign{\smallskip}\hline\noalign{\smallskip}
$H_{n}$ & $H_{n}(t) = \int_{V}|\nabla^{n}\bu|^{2}\,dV$ & $n\geq 0$\\
$\epsilon_{av}$ & $\epsilon_{av} = \nu L^{-3}\left<H_{1}\right>$ & \\
$\eta_{K}$ & $\eta_{K}^{-4} = \epsilon_{av}/\nu^{3}$ & \\
$F_{n}$ & $F_{n}(t) = \int_{V} \left(|\nabla^{n}\bu|^{2} + 
\tau^{2}|\nabla^{n}\bdf|^{2}\right)\,dV$ & \\
$\omega_{0}$ & $\omega_{0} = \nu L^{-2}$ & \\
$\tau$ & $\omega_{0}\tau = \Gr^{-(1/2 +\delta)}$ & $0 < \delta < \sixth$\\
$\kappa_{n}$ & $\kappa_{n}(t) = \left(F_{n}/F_{0}\right)^{1/2n}$ & \\
$\lambda_{n}$ & $\lambda_{n} = 3- \frac{5}{2n}+\frac{\delta}{n}$ & \\
\noalign{\smallskip}\hline
\end{tabular}
\ec
\end{table}

\section{Summary of the main results of the paper}\label{summary}

Sections 3-5 of this paper give a full account of the ideas with proofs from
first principles. This section has been introduced for readers who prefer to
peruse the proofs later. The next sub-section on long-time averages is
followed by a second describing where problems with regularity lie. The third
describes how the time-axis can be split into two types of interval, `good'
and `bad'. \S\ref{bad} summarises the dynamics on the bad intervals and the
very large lower bounds that can be found within these.

\subsection{Long-time averages}\label{longav1}

The best known long-time average in Navier-Stokes analysis is the explicit 
upper bound on the time averaged dissipation rate $\epsilon_{av} = \nu L^{-3}
\left<H_{1}\right>$ found from Leray's energy inequality \cite{Leray}
\bel{Ler1}
\shalf\dot{H}_{0} \leq -\nu H_{1} + H_{0}^{1/2}\|\bdf\|_{2}
\ee
Using Doering and Foias's result \cite{DF} that $\Gr \leq c\,\Rey^{1/2}$, 
it is easily shown that as $\Gr\to\infty$
\bel{Ler2}
\left<H_{1}\right> \leq c\,\nu^{2}L^{3}\ell^{-4}\Rey^{3}
\ee
Recall that the long-time average $\left<\cdot\right>$ is defined in (\ref{I8}).
The above $\Rey^{3}$ long-time average estimate for $\left<H_{1}\right>$ leads
to an upper bound on the energy dissipation rate $\epsilon_{av} \leq c\,\nu^{3}
\ell^{-4}\Rey^{3}$ and thence to the conventional estimate on the inverse 
Kolmogorov length $\eta_{K}^{-1}$ defined in (\ref{Koldef1}) 
\bel{Koldef2}
\eta_{K}^{-1}=\left(\frac{\epsilon_{av}}{\nu^3}\right)^{1/4}
\leq \ell^{-1}\Rey^{3/4}
\ee
Foias, Guillop\'{e} \& Temam's generalization of (\ref{Ler2}), when 
the inequality $\Gr\leq c\,\Rey^{1/2}$ is applied \cite{FGT}, can be 
expressed as (see Theorem \ref{FGTthm} in \S\ref{longav2}) 
\bel{FGT0}
\ell\left<F_{n}^{\frac{1}{2n-1}}\right> \leq c_{n}(L\ell^{-1})^{3}
\nu^{\frac{2n}{2n-1}}\Rey^{3}
\ee
When $n=1$, (\ref{FGT0}) recovers the sharp result of Doering and Foias \cite{DF} 
for $\left<F_{1}\right>$, except for a spurious volume factor $(L\ell^{-1})^{3}$ 
on the right hand side. 

A closely connected and important result used in this paper is the boundedness 
of the long-time averages of the $\kn$, for $n\geq 1$ taken with $\ell =L/2\pi$ 
(see \cite{Doering02} and Theorem \ref{knbd} in \S\ref{longav2})
\bel{I6}
\left<L\kn\right> \leq c_{n}\Rey^{\lambda_{n}}\hspace{2cm}
\lambda_{n}= 3-\frac{5}{2n} +\frac{\delta}{n}
\ee
where $\delta$ is a small parameter lying in the range $0<\delta <\sixth$.
From either ({\ref{FGT0}) or (\ref{I6}) an upper bound can also found on
$\left<\|\bu\|_{\infty}\right>$; see Theorem \ref{FGTthm} in \S\ref{longav2}
and Table \ref{tab5}. Note that the case $n=1$ gives $\left<L\kappa_{1}\right>
\leq c_{1}\Rey^{1/2}$, which is consistent with the traditional scaling of
$\Rey^{-1/2}$ for the Taylor micro-scale. Boundedness of the long-time average
of $\kn$ does not, of course, imply that the $\kn$ are bounded pointwise in
time.

\subsection{Problems with regularity: an illustration}\label{regprob}

Normal practice has been to consider the time evolution of the $F_{n}$ using 
differential inequalities.  With variations, this has been standard the approach 
taken since the early days of the subject \cite{Temam,CF,FMRT,DGbook}. One such 
inequality is
\bel{kladder0}
\shalf\dot{F}_{n} \leq -\shalf\nu F_{n+1} + 
\left(c_{n}\nu^{-1}\|\bu\|^{2}_{\infty} +  \nu \ell^{-2}\Rey\right)F_{n}
\ee
where $\|\bu(\cdot,\,t)\|_{\infty}$ is, in effect, the peak velocity on the
whole domain. The reader can find the precise derivation of (\ref{kladder0})
in Proposition \ref{Fnladder} in \S\ref{Fnsect}. The right hand side has two
dominant terms; one negative associated with the dissipation, and the dominant
positive $\|\bu\|_{\infty}^{2}$-term. Rewriting the $F_{n+1}$ term in
(\ref{kladder0}) in terms of the $\kn$ defined in (\ref{I4})
\bel{fnp1}
F_{n+1} = \kn^{2}\left(\frac{\kappa_{n+1}}{\kappa_{n}}\right)^{2(n+1)}F_{n}.
\ee
and using a Sobolev inequality
$\|\bu\|^{2}_{\infty}\leq c\,\kn^{3}F_{0}$ for $n\geq 2$, (\ref{kladder0}) becomes
\bel{kladder1}
\shalf\dot{F}_{n} \leq 
\left\{-\shalf\nu\kn^{2}\left(\frac{\kappa_{n+1}}{\kappa_{n}}\right)^{2(n+1)} + 
c_{n}\nu^{-1}\kn^{3}F_{0} +  \nu\ell^{-2}\Rey\right\}F_{n}.
\ee
From (\ref{I5}) the ratio $\knp/\kn$ has a lower bound of unity, $\knp/\kn
\geq 1$, for all periodic divergence-free functions. This reduces
(\ref{kladder1}) to
\bel{kladder2}
\shalf\dot{F}_{n} \leq \left(-\shalf\nu\kn^{2} + c_{n}\nu^{-1}\kn^{3}F_{0} +  
\nu\ell^{-2}\Rey\right)F_{n}.
\ee
The manifestly negative term $\sim\kn^{2}$ in (\ref{kladder2}) is not sufficient 
to control the $\kn^{3}$ term: arbitrarily large initial data on $F_{n}$ can be 
chosen that makes the right hand side positive. Despite the finiteness of the 
time averages (\ref{FGT0}) and (\ref{I6}), this leads to a failure to control 
either $F_{n}$ or $\kn$, other than for short times or small initial data 
\cite{Lady}.

\subsection{The potentially binary nature of the time-axis}\label{bin1}

One way of proceeding with this difficult problem is to ask whether the lower
bound of unity on $\knp/\kn$ could be raised, thereby effectively increasing
the dissipation in (\ref{kladder1}). Batchelor \& Townsend \cite{BT49}
experimentally identified a similar quantity to $\knp/\kn$ for $n=2,\,3$ that
was larger than expected for Gaussian data. The principal result of this
paper, proved in \S\ref{inter}, is that an improved lower bound on $\knp/\kn$
can be found which is valid only on sections of the time-axis. The estimates
in Theorem \ref{intervalthm} in \S\ref{inter}, the key result of the paper,
show that for Navier-Stokes weak solutions
\bel{I9}
\left<\left[c_{n}\left(\frac{\kappa_{n+1}}{\kappa_{n}}\right)\right]^{1/\mu -1}
-\left[(L\kappa_{n})^{\mu}\Rey^{-\lambda_{n}}\right]^{1/\mu -1}\right>\geq 0
\ee
where the real parameter $\mu$ can take any value in the range $0 < \mu < 1$.
The $c_{n}$ are the same as in (\ref{I6}). Given that the long-time average in 
(\ref{I9}) is non-negative means that there must be intervals of the time-axis, 
called \textit{good intervals}, where the inequality
\bel{I10}
c_{n}\left(\frac{\kappa_{n+1}}{\kappa_{n}}\right) 
\geq \left(L\kappa_{n}\right)^{\mu}\Rey^{-\lambda_{n}} 
\ee
holds. It is easily seen that when (\ref{I10}) is applied to (\ref{kladder1})
on these intervals the strength of the dissipation is increased. This applies at
small scales ($L\kappa_{n} > c_{n}\Rey^{\lambda_{n}/\mu}$) where the lower bound
on $\kappa_{n+1}/\kn$ in (\ref{I10}) is raised away from unity. The divisor
($F_{0}$) within $\kn$ is bounded both above and below. Thus (\ref{kladder1})
can be turned into a proper differential inequality in $F_{n}$; the reader can
turn to \S\ref{good} for details. The result, which is intuitively obvious
from (\ref{kladder1}), is that the negative dissipation term is stronger than 
the positive nonlinear term when
\bel{I10a}
\mu > \frac{1}{2(n+1)}
\ee
so no singularities are possible on these intervals (see \S\ref{good}). 

The integrand in (\ref{I9}) cannot be guaranteed to be always positive,
however; so intervals on the rest of the time-axis, where it could be
negative, are designated as \textit{bad intervals}. On these the reverse
inequality must be satisfied
\bel{I11}
c_{n}\left(\frac{\kappa_{n+1}}{\kappa_{n}}\right) 
< \left(L\kappa_{n}\right)^{\mu}\Rey^{-\lambda_{n}}.
\ee
Because its is always true that $\kappa_{n+1}\geq \kn$, then
\bel{I12}
L\kappa_{n}(t) > c_{n} \Rey^{\lambda_{n}/\mu}
\ee
on these intervals. It is, of course, possible that there are no bad intervals
and that the whole time-axis is `good'. The positivity of the average in 
(\ref{I9}), however, ensures that the complete time-axis cannot be `bad'. This 
paper is based on  the worst-case supposition that bad intervals exist and need 
to be dealt with accordingly. There appears to be no further information 
within (\ref{I9}) as to their distribution, which may depend on the value 
of $n$.

\subsection{Dynamics on the bad intervals}\label{bad}

\vspace{-1.5cm}
$$
\begin{minipage}[ht]{9cm}
\setlength{\unitlength}{.75cm}
\begin{picture}(11,11)(0,0)
\thicklines
\put(0,0){\vector(0,1){8}}
\put(0,0){\vector(1,0){10}}
\put(0,8.2){\makebox(0,0)[b]{$\kappa_{n}(t)$}}
\put(10.2,0){\makebox(0,0)[b]{$t$}}
\thinlines
\put(3,0){\line(0,1){8}}
\put(3.5,0){\line(0,1){8}}
\put(7,0){\line(0,1){8}}
\put(7.5,0){\line(0,1){8}}
\multiput(0,2)(.1,0){100}{.}
\multiput(0,6)(.1,0){30}{.}
\multiput(3.5,6)(.1,0){65}{.}
\put(-.75,2){\makebox(0,0)[b]{\scriptsize$Re^{\lambda_{n}}$}}
\put(5.3,2.1){\makebox(0,0)[b]{\scriptsize Long-time average}}
\put(5,6.3){\makebox(0,0)[b]{\scriptsize$L\kappa_{n}>Re^{\lambda_{n}/\mu}$}}
\put(3.9,6.4){\vector(-1,0){.72}}
\put(6.2,6.4){\vector(1,0){1.1}}
\put(-.75,6){\makebox(0,0)[b]{\scriptsize$Re^{\lambda_{n}/\mu}$}}
\put(3.25,-.70){\makebox(0,0)[b]{\scriptsize$(\Delta t)_{b}$}}
\put(5.25,-.70){\makebox(0,0)[b]{\scriptsize$(\Delta t)_{g}$}}
\put(5.75,-.25){\vector(1,0){1.25}}
\put(4.75,-.25){\vector(-1,0){1.2}}
\qbezier[250](0,1)(3,0)(3.01,8)
\qbezier(3.5,7.75)(5,-7)(7,7.75)
\qbezier[300](7.5,8)(7.5,0)(10,1)
\end{picture}
\end{minipage}
$$
\par\vspace{.5cm}\noindent
{\small \textbf{Figure 2:} A descriptive picture, not to scale, of 
good/bad intervals for some value of $n\geq 2$; constants have been 
omitted.}\label{knfig}
\par\vspace{0.5cm}\noindent
Various questions remain. Firstly, if the bad intervals exist, are they 
finite in width? The answer to this is in the affirmative. Estimates for the 
widths of bad intervals $\tb$ are displayed in Table \ref{tab5} (see Theorem
\ref{badwidththm} in \S\ref{widthest}) with upper bounds on $\mu$ lying in two
ranges together with a lower bound given in (\ref{I10}). It can be seen that
these widths are exceedingly small; in fact $\tb\to 0$ as $\Rey\to\infty$. The
large lower bound on $\kn$ in (\ref{I12}) indicates the predominance of high
wave-numbers within these very short intervals.

Figure \ref{knfig} is a descriptive picture of a typical distribution of 
good/bad 
intervals on the $t$-axis. It is drawn in such a way that $\kn(t)$ looks 
a relatively flat function. While this is probable, some artistic licence 
has been taken for the following reasons. The lower bound in (\ref{I12}) 
makes values of $\kn(t)$ much larger than the upper bound on the time 
average (\ref{I6}), so it has to spend relatively long intervals of time 
in the good to recompense. Nevertheless, this does not prove that it is 
quiescent in the central part of the interval. As is shown in \S\ref{inter}, 
there is enough freedom within the upper bound on $\kn$, under the 
constraint of the long-time average, to allow it to 
reach large enough values so that it can connect with the next bad interval. 
While we believe that it is likely that $\kn$ is flat in the central 
region, we have been unable to prove this. It is pathologically possible 
that enough fine structure might exist within the good interval without 
violating the average. 

The positions of the intervals may differ as $n$ varies so a new figure is 
needed for each $n$; Lemma \ref{Fmlemma} shows that if any \textit{one} $\kn$ 
is bounded then all are bounded. Only if the intersection of all bad intervals 
is itself `bad' is a singularity possible; \S\ref{badintersect} addresses this 
question. The set containing the intersection of all the bad intervals for all 
$n\geq 1$ is designated there as $\mathcal{S}^{(\infty)}$ and is related to 
Scheffer's potentially singular set \cite{Scheffer76a}. Taking this into 
account narrows the range of $\mu$ to $\thfth < \mu < \shalf$. 


The bad intervals (and their intersection) can be divided into sub-intervals; 
for the $i$th bad interval let us
take the $j$th sub-interval on which $\dot{F}_{n} \geq 0$. This we call the
$ij$th `dangerous' sub-interval of width $\tplus$. Singularities can occur on
any interval within these because $F_{n}$ is increasing, whereas sub-intervals
on which $\dot{F}_{n} \leq 0$ no singularities can occur because $F_{n}$ is
bounded ($n\geq 1$) by its initial value at the start of the sub-interval. 
Estimates for the width of dangerous sub-intervals, of width $\tplus$, are 
smaller than those for $\tb$ and can be found in Theorem \ref{subthm} of 
\S\ref{subint} (see Table \ref{tab5}). It is then possible to find very large 
\textit{lower} bounds within $\tplus$ on various quantities such as an 
equivalent of the inverse Kolmogorov length $\eta_{+}^{-1}$, based upon lower 
bounds on $F_{1}$ within $\tplus$; it is here where the intermittency of the
dissipation field shows up. Additionally, the peak velocity and the peak
velocity gradient matrix also have very large lower bounds. These are
displayed in Table \ref{tab3}.

Conclusions that can be drawn from Tables \ref{tab3}, \ref{tab4} 
and \ref{tab5} at the back of this section are:
\ben
\item The bad intervals are exceptionally narrow for large $\Rey$.

\item The action within these intervals is intense; the lower bound 
on $F_{1}$ within them illustrates the strength of the dissipation 
field there.
 
\item The dynamics within these intervals is so fine, even if no 
singularities occur, that they would be exceptionally difficult to 
resolve numerically. 
\een

Thirdly Theorem \ref{chebthm} bounds below, in the average sense, the 
ratio of the widths of the good and bad intervals
\bel{flat1} 
\frac{\overline{(\Delta t)}_{g}}{\overline{(\Delta t)}_{b}} 
\geq \Rey^{\lambda_{n}\left(\frac{1}{\mu}-1\right)}
\ee
showing that relative to the bad intervals, the widths of the good 
increase with $\Rey$.

Fourthly, the final result of \S\ref{condreg} displayed in Theorem 
\ref{condregthm} is a conditional regularity result. Assume that 
the energy $H_{0}(t)$ has a lower bound within the dangerous 
sub-intervals, $\Delta t_{+}^{(i,j)} = t - \tij$, of the form
\bel{conreg1}
H_{0}(t)\geq H_{0}(\tij)e^{-\omega_{0}\Rey\Delta t_{+}^{(i,j)}}.
\ee
Then solutions of the Navier-Stokes equations are regular there. $\tij$ 
is the initial time for the dangerous sub-interval. The very large initial 
conditions on the energy $H_{0}(\tij)$ at the junction of the intervals is 
the main obstacle to completing a regularity proof.
\begin{table}[ht]
\bc
\caption{\label{tab3}A comparison between the upper bounds on the long-time 
averages in column two (see Theorems \ref{knbd} and \ref{FGTthm} in 
\S\ref{longav2}) and the very large lower bounds in dangerous sub-intervals 
$\tplus$ in column three (see Theorem \ref{subthm} in \S\ref{subint}). 
Notice there the very large lower bounds on $F_{1}$ illustrating the 
strong intermittency in the dissipation field. All multiplicative constants 
have been omitted to save space. In addition $\ell = L/2\pi$.}
\begin{tabular}{lll} 
\hline\noalign{\smallskip}
\extrarowheight{9pt}
& Long-time average $\left<\cdot\right>$ & On $\tplus$ intervals \\
\noalign{\smallskip}\hline\noalign{\smallskip}
Energy moments & $\left<L\kappa_{n}\right>\leq \Rey^{\lambda_{n}}$ &
$L\kappa_{n}\geq \Rey^{4+b_{n}}$\\
Enstrophy & $L^{-3}\left<F_{1}\right>\leq \omega_{0}^{2}\Rey^{3}$ &
$L^{-3}F_{1}\geq \omega_{0}^{2}\Rey^{4+b_{n}}$\\
$\mbox{(Kolmogorov~scale)}^{-1}$ & $L\eta_{K}^{-1}\leq \Rey^{3/4}$ & 
$L\eta^{-1}_{+}\geq \Rey^{(4+b_{n})/4}$\\
$\mbox{(Taylor~micro-scale)}^{-1}$ & $\left<L\kappa_{1}\right>\leq \Rey^{1/2}$&
$L\kappa_{1}\geq \Rey^{b_{n}/2}$\\
Peak velocity & $\left<\|\bu\|_{\infty}\right>\leq L\omega_{0}\Rey^{3}$ &
$\|\bu\|_{\infty}\geq L\omega_{0}\Rey^{4+b_{n}}$\\
Vel. gradient matrix & $\left<\|\nabla\bu\|_{\infty}^{1/2}\right>^{2}\leq
\omega_{0}\Rey^{3}$ & $\|\nabla\bu\|_{\infty}\geq \omega_{0}\Rey^{4+b_{n}}$\\
\noalign{\smallskip}\hline
\end{tabular}
\ec
\end{table}
\begin{table}[ht]
\bc
\caption{\label{tab4}Definitions of $\lambda_{n},~ a_{n}$ and $b_{n}$. On the 
intersection set $\mathcal{S}^{(\infty)}$ of \S\ref{badintersect}, $\lambda_{n}$ 
can be replaced by $\Lambda^{(\infty)}_{n}$.}
\begin{tabular}{lll} 
\hline\noalign{\smallskip}
\extrarowheight{9pt}
 & Definition & Range\\
 \noalign{\smallskip}\hline\noalign{\smallskip}
$\lambda_{n}$ & $\lambda_{n} = 3- \frac{5}{2n}$ & \\
$a_{n}$ & $a_{n} = \frac{\lambda_{n+1}}{\mu}\left(\frac{2n-2}{2n-1}\right)-
\frac{10n-1}{2n-1}$ & $a_{n} >0$\\
$b_{n}$ & $b_{n} = \frac{\lambda_{n+1}}{\mu} - 4$ & $b_{n} >a_{n}$\\
\noalign{\smallskip}\hline
\end{tabular}
\ec
\end{table}
\begin{table}[ht]
\bc
\caption{\label{tab5}Estimates of widths of bad intervals, their 
sub-intervals $\tplus$, and corresponding ranges of $a_{n},~b_{n}$ 
\& $\mu$ (see Theorem \ref{badwidththm} in \S\ref{widthest} and 
Theorem \ref{subthm} in \S\ref{subint}). It is always true that 
$b_{n} > a_{n}$ (see Table \ref{tab4}). For the case $n=1$, $\mu 
> \shalf$.} 
\begin{tabular}{ccll}
\hline\noalign{\smallskip}
\extrarowheight{9pt}
Range of $\mu$ & $a_{n},~b_{n}$ & Width of interval $n\geq 2$ & \\
\noalign{\smallskip}\hline\noalign{\smallskip}
$\mu \leq \lambda_{n+1}\left(\frac{n-1}{6n-1}\right)$&$a_{n}\geq 1$
& $\omega_{0}(\Delta t)_{b} \leq \Rey^{-a_{n}}$ \\
$\lambda_{n+1}\left(\frac{n-1}{6n-1}\right) < \mu <
\lambda_{n+1}\left(\frac{2n-2}{10n-1}\right)$ & $ 0 < a_{n}< 1$ &
$\omega_{0}(\Delta t)_{b} \leq \Rey^{-1}\ln \Rey$\\
& $b_{n}> 1$ & $\omega_{0}\tplus \leq \Rey^{-b_{n}}$ &\\
$\frac{1}{2(n+1)} <\mu$ & & $\overline{(\Delta t)}_{g} \geq 
\overline{(\Delta t)}_{b}\, 
\Rey^{\lambda_{n}\left(\frac{1}{\mu}-1\right)}$& \\
\noalign{\smallskip}\hline
\end{tabular}
\ec
\end{table}

\section{Standard estimates}\label{standard}

\subsection{The forcing does not dominate the fluid}\label{forcing}

The technical parts of this paper revolve around the quantities 
\bel{Fndef}
F_{n}(t) = H_{n} + \tau^{2}\|\nabla^{n}\bdf\|^{2}_{2}
\ee
where the $H_{n}$ are defined in (\ref{Hndef}). The $F_{n}$ contain the fluid
velocity derivatives and those of the forcing, although the latter has been
assumed to have a narrow-band character, as shown in Table \ref{tab1}. They
are included within $F_{n}$ in order to circumvent problems that may arise
when dividing by these (squared) semi-norms. Once these terms have been 
introduced it is necessary to demonstrate that they do not dominate the 
fluid \cite{Doering02}, a point that has also made by Tsinober \cite{Tsinbk} 
in numerical computations. The 
characteristic time $\tau$ will be chosen for convenience but so long as 
$\tau \neq 0$, the $F_{n}$ are bounded away from zero by the explicit value 
$\tau^{2}L^{3}\ell^{-2n}f^2$. Moreover, $\tau$ may be chosen to depend on 
the parameters of the problem such that $\left<F_{n}\right>\sim \left<H_{n}\right>$ 
as $\Gr\to\infty$. To see how to achieve this, let us define
\bel{add1}
\tau = \ell^{2}\nu^{-1}\Gr^{-(\delta + 1/2)}
\ee
with $\delta > 0$, which is a parameter yet to be determined. Then the 
additional term in (\ref{Fndef}) is
\bel{add2}
\tau^{2}\|\nabla^{n}\bdf\|_{2}^{2} =
L^{3}\nu^{-2}\ell^{4-2n}f^{2}\Gr^{-(2\delta +1)}
= \nu^{2}\ell^{-(2n+2)}L^{3}\Gr^{1-2\delta}.
\ee
Now it has also been shown by Doering \& Foias \cite{DF} that the energy 
dissipation rate $\epsilon_{av} = \nu L^{-3}\left<H_{1}\right>$ has a 
\textit{lower} bound for high $\Gr$
\bel{Ler3}
\epsilon_{av} \geq c\,\nu^{3}\ell^{-3}L^{-1}\Gr
\ee
which can be used in the far right hand side of (\ref{add2}) 
\beq{add3}
\tau^{2}\|\nabla^{n}\bdf\|_{2}^{2} 
&\leq & c\,\epsilon_{av}\,\ell^{-(2n-1)}L^{4}\nu^{-1}\Gr^{-2\delta}\nonumber\\
& = & c\,\left(L\ell^{-1}\right)^{(2n-1)}
L^{-2(n-1)}\bigl<H_{1}\bigr>\Gr^{-2\delta}\nonumber\\
& \leq & c\,\left(L\ell^{-1}\right)^{(2n-1)}
\bigl<H_{n}\bigr>\Gr^{-2\delta}
\eq
where Poincar\'{e}'s inequality has been used at the last step. Hence, for 
any $\delta > 0$ the additional forcing term in (\ref{Fndef}) is seen to 
be negligible with respect to $\left<H_{n}\right>$ as $\Gr\to \infty$. The 
parameter $\delta$ is left arbitrary at this stage, although it will be 
restricted further in the course of proving results in the next section.
Our interest lies in results for high $\Gr$ so correction terms described 
above will be ignored and it can safely be said that for $\delta >0$ 
\bel{Ler4}
\left<F_{1}\right> \leq c\,\nu^{2}L^{3}\ell^{-4}\Rey^{3}
\ee
This is Leray's result for weak solutions with narrow-band forcing included 
in a rational manner; the next section deals with long-time averages of other 
quantities.

\subsection{The $F_{n}$-ladder}\label{Fnsect}

In the calculations that follow the $H_{n}$ are formally manipulated even
though they are not known to be finite pointwise in time for weak solutions;
the end results may be justified in the standard way by proceeding from a 
Galerkin approximation to the solutions and then removing the regularization 
in the final results. In the usual manner `$c$' and $c_{n}$ are used as 
generic constants.
\begin{proposition}\label{Fnladder}
For $\Gr\to\infty$ and $\delta >0$ the $F_{n}$ satisfy the 
differential inequalities
\bel{ladderlemma1}
\shalf\dot{F}_{0} \leq - \nu F_{1} + c_{1}\nu\ell^{-2}\Rey^{1+2\delta}F_{0},
\ee
\bel{ladderlemma2}
\shalf\dot{F}_{1} \leq - \squart\nu F_{2} + c_{2}\nu^{-3}F_{1}^{3} + 
c\,\nu\ell^{-2}\Rey^{1+2\delta}F_{1},
\ee
and, for $n\geq 2$, there are two alternative versions 
\bel{ladderlemma3}
\shalf\dot{F}_{n} \leq - \nu F_{n+1} + c_{n,1}\left(\|\nabla\bu\|_{\infty} 
+ \nu\ell^{-2}\Rey^{1+2\delta}\right)F_{n}.
\ee
\bel{ladderlemma4}
\shalf\dot{F}_{n} \leq - \shalf\nu F_{n+1} + c_{n,2}\left(\nu^{-1}
\|\bu\|_{\infty}^{2} + \nu\ell^{-2}\Rey^{1+2\delta}\right)F_{n}.
\ee
\end{proposition}
\begin{remark}
When the inequality $F_{1}^{2}\leq F_{2}F_{0}$ is used in
(\ref{ladderlemma2}), the resulting differential inequality for $F_{1}$
demonstrates the inability of these methods, as they stand, to gain control
over $F_{1}$ for arbitrarily large initial data.
\end{remark}
\begin{proof}
The proof follows in four steps. 
\par\medskip\noindent
\textbf{Step 1:}  Let us begin with the proof of (\ref{ladderlemma1}): Leray's 
energy inequality is
\bel{Enin1}
\shalf\dot{H}_{0}\leq -\nu H_{1} + H_{0}^{1/2}\|\bdf\|_{2}
\ee
Adding and subtracting the quantity $\nu\tau^{2}\|\nabla\bdf\|_{2}^{2}$, it is 
seen that
\bel{Enin2}
\shalf\dot{F}_{0}\leq -\nu F_{1} + \nu\tau^{2}\|\nabla\bdf\|_{2}^{2} + 
H_{0}^{1/2}\|\bdf\|_{2}
\ee
Because the forcing is narrow-band as in Table \ref{tab1}, it is possible to
reduce a derivative on the $\|\nabla\bdf\|_{2}^{2}$ term. This, together with
Young's inequality (using $g\tau^2 > 0$ as a parameter with $g$ is to be
suitably chosen below) to break up the last term
\bel{Enin3}
\shalf\dot{F}_{0}\leq -\nu F_{1} + \frac{1}{2g\tau^2}H_{0} + 
\tau^{2}\left(\shalf g + \frac{\nu}{\ell^{2}}\right)\|\bdf\|_{2}^{2}
\ee
$g$ is determined by making the coefficients of $H_{0}$ and
$\tau^{2}\|\bdf\|_{2}^{2}$ equal
\bel{Enin4}
g = - \frac{\nu}{\ell^2} + \left\{\nu^{2}\ell^{-4} + \tau^{-2}\right\}^{1/2}.
\ee
With $\tau$ chosen as in (\ref{add1}) with $\delta > 0$, $g$ becomes
\bel{Enin5}
g = \tau^{-1}\left(\left\{1+\Gr^{-(2\delta +1)}\right\}^{1/2} - 
\Gr^{-(\delta + 1/2)}\right).
\ee
Consequently, $g\sim \tau^{-1}$ as $\Gr\to\infty$ in which case  
\bel{taudef}
\tau^{-1} = \nu\ell^{-2}\Gr^{\shalf +\delta}\leq c\,\nu\ell^{-2}\Rey^{1+2\delta}
\ee
In this limit (\ref{Enin3}) can be written as in (\ref{ladderlemma1}). 
\par\medskip\noindent
\textbf{Step 2:} The proof of (\ref{ladderlemma2}) is found directly 
from ($\bw = \mbox{curl}\,\bu$)
\bel{F1a}
\shalf\dot{H}_{1} \leq - \nu H_{2} + 
\int_{\Omega}\bw\cdot(\bw\cdot\nabla)\bu\,dV 
+ H_{1}^{1/2}\|\nabla\bdf\|_{2}.
\ee
The middle term can be estimated thus
\bel{F1b}
\int_{\Omega}\bw\cdot(\bw\cdot\nabla)\bu\,dV \leq
\|\bw\|_{4}^{2}\|\nabla\bu\|_{2}\leq c\,F_{2}^{3/4}F_{1}^{3/4}
\ee
having used the Sobolev inequality $\|\bw\|_{4}\leq c\,\|\nabla\bw\|_{2}^{3/4}
\|\bw\|_{2}^{1/4}$. The procedure with the forcing is then used as in Step 1 to 
obtain
\bel{ex1}
\shalf\dot{F}_{1} \leq - \nu F_{2} + c\,F_{2}^{3/4}F_{1}^{3/4} + 
c\,\nu\ell^{-2}\Rey^{1+2\delta}F_{1},
\ee
Young's inequality on the middle term finally gives (\ref{ladderlemma2}).
\par\medskip\noindent
\textbf{Step 3:} For a proof of (\ref{ladderlemma3}), consider the ladder of
differential inequalities satisfied by the $H_{n}$ for $n\geq 2$ (see
\cite{DGbook,BDG})
\bel{ladder2}
\shalf\dot{H}_{n} \leq - \nu H_{n+1} + c_{n}\|\nabla\bu\|_{\infty}H_{n} + 
H_{n}^{1/2}\|\nabla^{n}\bdf\|_{2}.
\ee
(\ref{ladderlemma3}) is proved by following the procedure with the forcing 
as in Step 1. 
\par\medskip\noindent
\textbf{Step 4:}  The alternative to the differential inequality (\ref{ladder2})
for $n\geq 2$ is
\bel{ladder3}
\shalf\dot{H}_{n} \leq - \shalf\nu H_{n+1} + 
c_{n}\nu^{-1}\|\bu\|_{\infty}^{2}H_{n} 
+ H_{n}^{1/2}\|\nabla^{n}\bdf\|_{2}.
\ee
(\ref{ladderlemma4}) is found by using the same procedure as in Step 1 except
that the quantity $\shalf\nu\tau^{2}\|\nabla^{n+1}\bdf\|_{2}^{2}$ is 
subtracted whereas $\nu\tau^{2}\|\nabla^{n+1}\bdf\|_{2}^{2}$ is 
added.
\end{proof}
\subsection{Long-time averages}\label{longav2}

In \S\ref{intro} it was shown how the quantities $\kappa_{n}$ are ordered
such that $\kn \leq \kappa_{n+1}$. There is no $\kappa_{n}$ that is known to
be \textit{a priori} bounded. What is known is the boundedness of the long-time 
averages defined in (\ref{I8}). The equivalent of Leray's bulk dissipation 
estimate in terms of $\kappa_{1}$ instead of $F_{1}$ is found from 
(\ref{ladderlemma1}) by dividing through $F_{0}$ and long-time averaging
\bel{Ler5}
\ell^{2}\left<\kappa_{1}^{2}\right> \leq c\,\Rey^{1+2\delta}.
\ee
The first of the following two theorems states results on long-time averages for 
higher values of $n$. This estimate can be found in \cite{Doering02} with a 
wider range of $\delta$.
\begin{theorem}\label{knbd}
For $\Gr\to\infty$ and the parameter $\delta$ lying in the range $0 < \delta < 
\sixth$
\bel{kappalemma2}
\ell\left<\kn\right> \leq c_{n}\left(L\ell^{-1}\right)^{\frac{3(n-1)}{n}}
\Rey^{\lambda_{n}} \hspace{2cm}n\geq 1
\ee
where $\lambda_{n}$ is defined by 
\bel{lamdef}
\lambda_{n} = 3-\frac{5}{2n} + \frac{\delta}{n}.
\ee
\end{theorem}
\begin{proof}
\textbf{Step 1:} Consider first $\kappa_{2,1}$:
\beq{k1}
\left<\kappa_{2,1}\right> &=& \left<\left(\frac{F_{2}}{F_{1}}\right)^{1/2}\right>
\leq \left<\frac{F_{2}}{F_{1}^2}\right>^{1/2}\left<F_{1}\right>^{1/2}\nonumber\\
&\leq& \frac{\nu^2}{2}\left<\frac{F_{2}}{F_{1}^2}\right> 
+ \frac{1}{2\nu^{2}}\left<F_{1}\right>
\eq
where Young's inequality has been used at the last step.  Dividing inequality 
(\ref{ladderlemma2}) in Proposition \ref{Fnladder} in \S\ref{Fnsect} by 
$F_{1}^{2}$ and long-time averaging gives 
\bel{k2}
\nu^{2}\left<\frac{F_{2}}{F_{1}^2}\right> \leq 
c\,\nu^{-2}\left<F_{1}\right> + \nu\tau^{-1}\left<F_{1}^{-1}\right>
\ee
and so
\bel{k3}
\left<\kappa_{2,1}\right> 
\leq c\,\nu^{-2}\left<F_{1}\right> + \nu\tau^{-1}\left<F_{1}^{-1}\right>
\ee
The last term is
\bel{k4}
\nu\tau^{-1}\left<F_{1}^{-1}\right> 
\leq \frac{\nu\tau^{-1}}{\tau^{2}\ell^{-2}L^{3}f^{2}}
= \ell^{2}L^{-3}\Gr^{3\delta -\shalf}
\ee
from which it is concluded that $\delta$ must lie in the range $0 < \delta <
\sixth$ to be certain that this term decreases as $\Gr\to\infty$. Because
\bel{k5}
\left<F_{1}\right> \leq \nu^{2}L^{3}\ell^{-4}\Rey^{3}
\ee
then it follows that
\bel{k6}
\ell\left<\kappa_{2,1}\right> \leq c\,\left(L\ell^{-1}\right)^{3}\Rey^{3}
\ee
\par\medskip\noindent
\textbf{Step 2:} Now consider the quantities $\left<\kappa_{n+1,n}\right>$ for
$n\geq2$
\beq{k7}
\left<\kappa_{n+1,n}\right> &=& 
\left<\left(\frac{F_{n+1}}{F_{n}^{2n/(2n-1)}}\right)^{1/2}F_{n}^{1/2(2n-1)}\right>
\nonumber\\
&\leq& \left<\frac{F_{n+1}}{F_{n}^{2n/(2n-1)}}\right>^{1/2}
\left<F_{n}^{1/(2n-1)}\right>^{1/2}
\eq
\beq{k8}
\left<F_{n}^{\frac{1}{2n-1}}\right> &=& \left<\kappa_{n,1}^{(2n-2)/(2n-1)}
F_{1}^{1/(2n-1)}\right>\nonumber\\
&\leq& \left<\kappa_{n,1}\right>^{(2n-2)/(2n-1)}\left<F_{1}\right>^{1/(2n-1)}
\eq
Having used the fact that $\kappa_{n,1}\leq \kappa_{n+1,n}$, (\ref{k7}) and 
(\ref{k8}) give 
\bel{k8a}
\left<\kappa_{n+1,n}\right> \leq \left[\nu^{\frac{2}{2n-1}}
\left<\frac{F_{n+1}}{F_{n}^{(2n-1)/2n}}\right>\right]^{\frac{(2n-1)}{2n}}
\left[\nu^{-2}\left<F_{1}\right>\right]^{\frac{1}{2n}}
\ee
so a H\"{o}lder inequality gives 
\bel{k9}
2n\left<\kappa_{n+1,n}\right> \leq (2n-1)\nu^{\frac{2}{2n-1}}
\left<\frac{F_{n+1}}{F_{n}^{2n/(2n-1)}}\right> + \nu^{-2}\left<F_{1}\right>
\ee
To estimate the first long-time average on the right hand side, consider the 
second $F_{n}$-ladder in (\ref{ladderlemma4}) 
\bel{k10}
\shalf\dot{F}_{n} \leq - \shalf\nu F_{n+1} + c_{n}\left(\nu^{-1}
\|\bu\|_{\infty}^{2} + \nu\ell^{-2}\Rey\right)F_{n}.
\ee
Now define
\bel{Yndef}
Y_{n} = F_{n}^{-\frac{1}{2n-1}}
\ee
and turn (\ref{k10}) into a differential inequality in $Y_{n}$ which involves 
dividing by $F_{n}^{2n/(2n-1)}$. To achieve this we use $\|\bu\|_{\infty}^{2} 
\leq c\,\kappa_{2,1}F_{1}$ and recall that $\kappa_{2,1}\leq\kappa_{n,1}$, then
\bel{k10a}
\|\bu\|_{\infty}^{2} F_{n}^{-\frac{1}{2n-1}} \leq
c\,\kappa_{2,1}\left[\kappa_{n,1}^{-1}F_{1}\right]^{\frac{2n-2}{2n-1}}
\leq c\,\kappa_{2,1}^{\frac{1}{2n-1}}F_{1}^{\frac{2n-2}{2n-1}}
\ee
Hence (\ref{k10}) can be rewritten as 
\bel{k11}
(n-\shalf)(\dot{Y}_{n} + \nu\ell^{-2}\Rey Y_{n}) \geq 
\shalf\nu\frac{F_{n+1}}{F_{n}^{\frac{2n}{2n-1}}} - c\,\nu^{-1}
\kappa_{2,1}^{\frac{1}{2n-1}}F_{1}^{\frac{2n-2}{2n-1}}
\ee
Making up the coefficient in $\nu$ to that in (\ref{k9}), a H\"{o}lder
inequality on the last term gives
\beq{k12}
\nu^{\frac{2}{2n-1}} \frac{F_{n+1}}{F_{n}^{2n/(2n-1)}}
&\leq & (2n-1)\nu^{\frac{3-2n}{2n-1}}\left[\dot{Y}_{n} + 
\nu\ell^{-2}\Rey Y_{n}\right]\nonumber\\
&+& \frac{1}{2n-1}\left\{\kappa_{2,1} + c\,(2n-2)\nu^{-2}F_{1}\right\}.
\eq
Taking the long-time average of this in (\ref{k9}) we have 
\bel{k13}
2n\left<\kappa_{n+1,n}\right> \leq \left<\kappa_{2,1}\right> + 
c_{n}\,\nu^{-2}\left<F_{1}\right> + 
\nu^{\frac{2}{2n-1}}(2n-1)\ell^{-2}\Rey\left<Y_{n}\right>
\ee
The long-time average of $\dot{Y}_{n}$ has vanished and the last term 
$\left<Y_{n}\right>$ is bounded above (because $F_{n}$ is bounded below) so the 
long-time average is zero. Thus when (\ref{k5}) and (\ref{k6}) are used we have
\bel{kappalemma1}
\ell\left<\kappa_{n,1}\right>\leq
\ell\left<\kappa_{n+1,n}\right> \leq c_{n}\left(L\ell^{-1}\right)^{3}\Rey^{3}
\ee
Note that the exponents of $\Rey$ and $L/\ell$ are uniform in $n$; only the
constant is not. (\ref{kappalemma1}) can now be used to estimate
$\left<\kappa_{n}\right>$ in the final step.
\par\medskip\noindent
\textbf{Step 3:} Rewrite $\left<\kappa_{n}\right>$ in the following way: 
\beq{k14}
\left<\kappa_{n}^{\frac{2n}{2n-1}}\right> & = &
\left<\left(\frac{F_{n}}{F_{0}}\right)^{\frac{1}{2n-1}}\right>
= \left<\left(\frac{F_{n}}{F_{1}}\right)^{\frac{1}{2n-1}}
(\kappa_{1}^{2})^{\frac{1}{2n-1}}\right>\nonumber\\
&=& \left<\kappa_{n,1}^{\frac{2n-2}{2n-1}}
(\kappa_{1}^{2})^{\frac{1}{2n-1}}\right> \leq 
\left<\kappa_{n,1}\right>^{\frac{2n-2}{2n-1}}
\left<\kappa_{1}^{2}\right>^{\frac{1}{2n-1}}
\eq
Using our estimate for $\left<\kappa_{n,1}\right>$ from (\ref{kappalemma1}) 
and also that for $\left<\kappa_{1}^{2}\right>$ from (\ref{Ler5}), the result 
in (\ref{kappalemma2}) is proved.
\end{proof}
The first infinite set of non-trivial, bounded, long-time averages were 
those found by Foias, Guillop\'{e} \& Temam \cite{FGT}. These are related to 
those in 
Theorem \ref{knbd}, and particularly to the estimates for $\kappa_{n,1}$ in
(\ref{kappalemma1}).
\begin{theorem}\label{FGTthm}
For $\Gr\to\infty$ the long-time averaged quantities of Foias, Guillop\'{e} \& 
Temam \cite{FGT} are estimated in terms of $\Rey$ as
\bel{FGT1}
\ell\left<\|\bu\|_{\infty}\right> 
\leq c_{1}\nu\left(L\ell^{-1}\right)^{3}\Rey^{3} 
\ee
\bel{FGT2}
\ell\left<F_{n}^{\frac{1}{2n-1}}\right> 
\leq c_{n,2}\nu^{\frac{2}{2n-1}}\left(L\ell^{-1}\right)^{3}\Rey^{3} 
\ee 
\bel{FGT3}
\ell\left<\|\nabla\bu\|_{\infty}^{1/2}\right> 
\leq c_{3}\nu^{1/2}\left(L\ell^{-1}\right)^{3}\Rey^{3}.
\ee
\end{theorem}
\par\medskip\noindent
\begin{remark}
The case $n=1$ is distinct from the result in \cite{DF} because of 
the $\left(L\ell^{-1}\right)^{3}$ on the right hand side.
\end{remark}
\begin{proof}
The proof follows from the Sobolev inequalities 
\bel{FGT4}
\|\bu\|_{\infty} \leq c\,\kappa_{n,1}^{1/2}F_{1}
\hspace{2cm}
\|\nabla\bu\|_{\infty} \leq c\,\kappa_{n,1}^{3/2}F_{1}^{1/2}
\ee
with the estimates (\ref{kappalemma1}) for $\kappa_{n,1}$ and (\ref{Ler4}) for
$\left<F_{1}\right>$. The quantities in (\ref{FGT2}) can be rewritten in terms
of $\kappa_{n,1}$ and $F_{1}$ and the result follows.
\end{proof}
A lemma is now proved that will be useful in later sections:
\begin{lemma}\label{Fmlemma}
If any $F_{m}$ ($\kappa_{m}$) is bounded on a time interval $[0,\,T]$ for 
$1 \leq m \leq n$ then so are all $F_{n}$ ($\kn$) for $n>m$.
\end{lemma}
\begin{proof}
Consider (\ref{ladderlemma3}) in Proposition \ref{Fnladder} 
in \S\ref{Fnsect} above; for $n\geq 3$ a Sobolev inequality gives 
\bel{lem1}
\|\nabla\bu\|_{\infty} \leq c\, \|\nabla^{n}\bu\|_{2}^{a}\|\nabla\bu\|_{2}^{1-a}
\leq F_{n}^{a/2}F_{1}^{(1-a)/2}
\ee
where $a= 3/[2(n-1)]$. There is an inequality for the $F_{n}$ of the form 
\bel{lem2}
F_{N}^{p+q}\leq F_{N-p}\,F_{N+q}^{p}
\ee
The choice of $N = n,~p =n-1$ and $q=1$ gives
\bel{lem3}
-F_{n+1}\leq -F_{n}^{\frac{n}{n-1}}/F_{1}^{\frac{1}{n-1}}
\ee
so, in consequence, (\ref{ladderlemma3}) becomes
\bel{lem4}
\shalf\dot{F}_{n}\leq -\nu F_{n}^{\frac{n}{n-1}}/F_{1}^{\frac{1}{n-1}}
+c\,F_{n}^{1+ a/2}F_{1}^{(1-a)/2} + c\,\nu\ell^{-2}F_{n}
\ee
Because $n/(n-1) > 1+ a/2$, (\ref{lem4}) makes it clear that if $F_{1}$ is 
bounded above at any time then all $F_{n}$ are bounded. If any $F_{m}$ 
is bounded for $m>1$ then $F_{1}$ must also be bounded (from (\ref{lem2})), 
in which case all $F_{n}$ are bounded for any $n>m$. The same results hold 
for the $\kn$ because the divisor $F_{0}$ is bounded from above and below. 
\end{proof}

\section{Intermittency: the binary form of the time-axis}\label{inter}

In the summary section, \S\ref{summary}, it was discussed how the effective
viscosity could be increased by proving that the ratio $\kappa_{n+1}/\kn$ has
a lower bound that is greater than unity under certain circumstances. This was
discussed in the context of the ladder of differential inequalities
(\ref{kladder1}) for the $F_{n}$ which is repeated here
\bel{Fninequal}
\shalf\dot{F}_{n} \leq 
\left(-\shalf\nu\kn^{2}\left(\frac{\kappa_{n+1}}{\kappa_{n}}\right)^{2(n+1)} + 
c_{n}\nu^{-1}\kn^{3}F_{0} + \nu \ell^{-2}\Rey\right)F_{n}.
\ee
The task of this section is to investigate lower bounds on the ratio
$\kappa_{n+1}/\kappa_{n}$. In the rest of this paper, the two lengths $L$ 
and $\ell$ will be taken such that  $\ell = L/2\pi$ to reduce algebra. 
Additionally, the parameter $\delta$, lying in the range $0 < \delta 
< \sixth$, that appears in the exponents of many of the estimates of 
the previous section, will be taken as arbitrarily small (but fixed) 
and ignored hereafter. 
\begin{theorem}\label{intervalthm}
For the parameter $\mu$ taking any value in the range $0 <\mu <1$, the 
ratio $\kappa_{n+1}/\kn$ obeys the long-time averaged inequality $(n\geq 1)$
\bel{as4}
\left<\left[c_{n}\left(\frac{\kappa_{n+1}}{\kappa_{n}}\right)\right]^{1/\mu -1}
-\left[(L\kappa_{n})^{\mu}\Rey^{-\lambda_{n}}\right]^{1/\mu -1}\right>\geq 0
\ee
where the $c_{n}$ are the same as those in (\ref{kappalemma2}).  Hence there 
exists at least one interval of time, designated as a `good interval', on which 
the inequality
\bel{thm1}
c_{n}\left(\frac{\kappa_{n+1}}{\kappa_{n}}\right) 
\geq \left(L\kappa_{n}\right)^{\mu}\Rey^{-\lambda_{n}} 
\ee
holds. Those other parts of the time-axis on which the reverse inequality 
holds
\bel{thm2}
c_{n}\left(\frac{\kappa_{n+1}}{\kappa_{n}}\right) < 
\left(L\kappa_{n}\right)^{\mu}\Rey^{-\lambda_{n}}
\ee
are designated as 'bad intervals'.
\end{theorem}
\begin{remark}
In principle, the whole time axis could be a good interval, whereas the 
positive time average in (\ref{as4}) ensures that the complete time-axis 
cannot be `bad'. This paper is based on the worst-case supposition that 
bad intervals exist, that they could be multiple in number, and that the 
good and the bad are interspersed. This is what is meant in this paper by 
a `potentially binary character', although the precise distribution and 
occurrence of the good/bad intervals and how they depend on $n$ remains an 
open question. It will be left until later (Theorem \ref{badwidththm}) to 
prove that the bad intervals are finite in width.
\end{remark}
\begin{proof}
Take two parameters $0 < \mu < 1$ and $0 < \alpha < 1$ such 
that $\mu + \alpha = 1$. The inverses $\mu^{-1}$ and $\alpha^{-1}$ will be 
used as exponents in the H\"{o}lder inequality on the far right hand side of
\bel{as2}
\left<\kappa_{n}^{\alpha}\right> \leq
\left<\kappa_{n+1}^{\alpha}\right> = \left<\left(\frac{\kappa_{n+1}}
{\kappa_{n}}\right)^{\alpha}\kappa_{n}^{\alpha}\right> 
\leq 
\left<\left(\frac{\kappa_{n+1}}{\kappa_{n}}\right)^{\alpha/\mu}
\right>^{\mu}\left<\kappa_{n}\right>^{\alpha}
\ee
thereby giving 
\bel{as3}
\left<\left(\frac{\kappa_{n+1}}{\kappa_{n}}\right)^{\alpha/\mu}\right>
\geq \left(\frac{\left<\kappa_{n}^{\alpha}\right>}
{\left<\kappa_{n}\right>^{\alpha}}\right)^{1/\mu} 
= \left<\kappa_{n}^{\alpha}\right> \left(\frac{\left<\kappa_{n}^{\alpha}\right>}
{\left<\kappa_{n}\right>}\right)^{\alpha/\mu}.
\ee
Navier-Stokes information can be injected into these formal manipulations: the
weak solution upper bound (\ref{kappalemma2}) and the lower bound $L\kappa_{n}
\geq 1$ can be used in the ratio on the far right hand side of (\ref{as3}) to
give (\ref{as4}), with the same $c_{n}$ as in
(\ref{kappalemma2}).
\end{proof}

\subsection{Bounds within good intervals}\label{good}

On the good intervals, application of the improved lower bound (\ref{thm1}) to
the differential inequality (\ref{Fninequal}) appears to imply that $\mu$ must
satisfy $2\mu(n+1) >1$ for the exponent of the negative term to be larger than
the positive. To harden this argument it is necessary to convert
(\ref{Fninequal}) into a differential inequality in $F_{n}$ alone. This can be
achieved because the divisor within $\kn$, namely $F_{0}$, is bounded above
\beq{good1}
\shalf \dot{F}_{n} &\leq& -\nu \Rey^{-2\lambda_{n}(n+1)}
L^{2\mu(n+1)}F_{n}^{\frac{(1+\mu)(n+1)}{n}}
F_{0}^{-\frac{\mu(n+1)+1}{n}}\nonumber\\
& + & c\,\nu^{-1}F_{n}^{\frac{2n+3}{2n}}
F_{0}^{\frac{2n-3}{2n}} + \nu L^{-2}\Rey F_{n}
\eq
For $n\geq 2$ arbitrarily large initial data a singularity can be 
prevented from forming if the exponent of the negative $F_{n}$-term 
be greater than that of the positive
\bel{good2}
\frac{(1+\mu)(n+1)}{n} > \frac{2n+3}{2n}\hspace{1cm}
\Rightarrow \hspace{1cm}\mu > \frac{1}{2(n+1)}
\ee
as predicted. It is not possible to take the infinite time limit because of the 
finiteness of the interval but the value of $F_{n} = F_{n,max}$ that turns the 
sign of the right hand side of (\ref{good1}) is bounded by
\bel{Fnmax}
F_{n,max} \leq L^{-2n}\Rey^{\gamma_{n}}F_{0,max}\equiv U_{bd}
\ee
\bel{gamdef}
\gamma_{n} = \frac{4n[\lambda_{n}(n+1)+2]}{2\mu(n+1)-1}
\ee
The exponent $\gamma_{n}>0 $ when $2\mu (n+1) >1$ and we have used the fact 
that $F_{0,max} = c\,L\nu^{2}\Rey^{4}$. In terms of Figure \ref{knfig} 
in \S\ref{summary}, it is necessary to prove that the solution in
the good region can become large enough to form an initial condition for weak
solutions in the bad region. This can be proved by the following argument:
consider that on bad intervals the $\kappa_{n}$ are bounded below uniformly by
\bel{below}
\left[L\kappa_{n}(t)\right]^{\mu} > c_{n}\Rey^{\lambda_{n}}
\ee
where $c_{n,\mu} = c_{n}^{1/\mu}$. In terms of $F_{n}$ this can be 
expressed as
\bel{Fnmaxbelow}
F_{n} > c_{n,\mu}L^{-2n}\Rey^{\frac{2n\lambda_{n}}{\mu}}F_{0,min}
\equiv L_{bd}
\ee
The question revolves around the relative sizes of the lower bound $L_{bd}$ 
in (\ref{Fnmaxbelow}) and $U_{bd}$ in (\ref{Fnmax})
\bel{ratio}
\frac{U_{bd}}{L_{bd}}
= \left(\frac{F_{0,max}}{F_{0,min}}\right)
\Rey^{\frac{2n(\lambda_{n}+4\mu)}{\mu[2\mu(n+1)-1]}} > 1
\hspace{2cm}\Rey\gg 1
\ee
Hence it is possible for $F_{n}$ to reach magnitudes at the edges of the good 
region that lie above the lower bound in (\ref{Fnmaxbelow}). 

For the case $n=1$, the following Lemma is applicable
\begin{lemma}\label{n=1lemma}
When $n=1$ no singularity can form on good intervals provided $\mu>\shalf$.
\end{lemma}
\begin{proof}
This follows immediately by applying Theorem \ref{intervalthm} to 
(\ref{ladderlemma2}).
\end{proof}
Nothing has yet been proved so far regarding the widths of the good and bad
intervals, $\tb$ and $\tg$ respectively, nor have we any further information 
regarding their nature. While it is possible that they may form pathological 
fractal 
subsets of the time-axis it will be assumed that these intervals are simple 
open or closed sets; the next section is devoted to estimating upper bounds 
on $\tb$. Here it is shown that a \textit{lower} bound can be found on the 
ratio of the average widths of the good and bad intervals. The argument is 
based on an elementary application of the Markov-Chebychev inequality. 
Consider an interval of time $[t_{p},\,t_{q}]$ that contains an equal 
number $N$ of good and bad intervals of widths $\tg^{(i)}$ and $\tb^{(i)}$ 
respectively. Define the average widths as
\bel{cheb1}
\overline{(\Delta t)}_{g} = \lim_{N\to\infty}\frac{1}{N}\sum^{N}_{i=1}\tg^{(i)}
\hspace{1cm}
\overline{(\Delta t)}_{b} = \lim_{N\to\infty}\frac{1}{N}\sum^{N}_{i=1}\tb^{(i)}
\ee
\begin{theorem}\label{chebthm}
Consider an interval of time $[t_{p},\,t_{q}]$ containing $N$ pairs of good 
and bad intervals. In the limits $N\to\infty$ and $[t_{p},\,t_{q}]\to\infty$, 
provided $(\Delta t)_{b}>0$, the ratio $\overline{(\Delta t)}_{g}/
\overline{(\Delta t)}_{b}$ diverges as
\bel{cheb2}
\frac{\overline{(\Delta t)}_{g}}{\overline{(\Delta t)}_{b}}
\geq c_{n}\Rey^{\lambda_{n}\left(\frac{1}{\mu}-1\right)}
~~~~~~\mbox{as}~~~\Rey\to\infty
\ee
\end{theorem}
\begin{proof}
Given (\ref{cheb1}), the fraction of time occupied by the bad 
intervals satisfies
\beq{cheb3}
\frac{\sum^{N}_{i=1}\tb^{(i)}}{\sum^{N}_{i=1}[\tg^{(i)}+\tb^{(i)}]} &\leq
& \frac{1}{t_{q}-t_{p}}\int_{T_{p,q}}dt \nonumber\\
&\leq & \frac{1}{t_{q}-t_{p}}\left(\frac{\int_{[t_{p},\,t_{q}]}L\kappa_{n}dt}
{c_{n,\mu}\Rey^{\lambda_{n}/\mu}}\right).
\eq
where $T_{p,q} = [L\kappa_{n}(t)\geq c_{n,\mu}\Rey^{\lambda_{n}/\mu}]
\cap [t_{p},\,t_{q}]$, so as $N\to\infty$ and $t_{q}-t_{p}\to\infty$, we have
\bel{cheb4}
\frac{\overline{(\Delta t)}_{b}}{\overline{(\Delta t)}_{g}+\overline{(\Delta
t)}_{b}} \leq \frac{\left<L\kappa_{n}\right>} {c_{n,\mu}\Rey^{\lambda_{n}/\mu}}
\leq \left[c_{n}\Rey^{\lambda_{n}}\right]^{1-\frac{1}{\mu}}
\ee
where we have used (\ref{below}) and (\ref{kappalemma2}). Hence we have the 
result.
\end{proof}

\section{What happens in the bad intervals?}\label{badinter}

It is necessary to prove that the bad intervals are of finite width: that is,
an upper bound is required on $\Delta t = t - t_{0}$ where $t_{0}$ is the
initial time of some arbitrary bad interval. Technically speaking, there
should be a superscript label for the $i$th bad interval such that $\Delta t
\equiv \Delta t^{(i)}$ and another on $t_{0}\equiv t_{0}^{(i)}$, but these
have been dropped for convenience. Recall that $\omega_{0} = \nu L^{-2}$ and
\bel{Edef1}
\calE (\Delta t ) = \frac{e^{\omega_{0}\Rey\,\Delta t}-1}{\omega_{0}\Rey} 
\ee
It will become necessary to solve inequalities of the type 
\bel{Einequal1}
\cE(\Delta t ) \leq \omega_{0}^{-1}\Rey^{-\beta}
\ee
for $\beta >0$ as $\Rey\to\infty$. It is not difficult to show that when 
$\beta \geq 1$, to leading order 
\bel{Einequal2}
\omega_{0}(\Delta t ) \lesssim \Rey^{-\beta}
\ee
whereas when $0 < \beta < 1$ then, to leading order
\bel{Einequal3}
\omega_{0}(\Delta t ) \lesssim (1-\beta)\Rey^{-1}\ln \Rey
\ee
The main task of this section is to show that the bad intervals have a finite 
widths and to find an upper bound on these. This requires two subsidiary 
estimates for
\bel{Edef2}
\Idt e^{\omega_{0}\Rey\,\Delta t}F_{1}(t)\,dt
\hspace{1cm}\mbox{and}\hspace{1cm}
\Idt e^{\omega_{0}\Rey\,\Delta t}\kappa_{2,1}(t)\,dt
\ee

\subsection{Two subsidiary estimates}\label{subsid}

\begin{lemma}\label{lemma3}
An estimate for the exponentially weighted time integral of $F_{1}$ is
\bel{F1est1}
\Idt e^{\omega_{0}\Rey\,\Delta t}F_{1}(t)\,dt \leq  
c_{1}\nu L \Rey^{4} + c_{2}\nu^{2}L^{-1}\calE(\Delta t)
\left[\Rey^{5}+O(\Rey^{4})\right]
\ee
\end{lemma}
\begin{proof}
Let us denote $H_{0}(t) = X^{2}(t)$ with $H_{0}(t_{0}) =
X^{2}_{0}$ then Leray's energy inequality (\ref{Ler1}) for weak solutions, 
\bel{Ler1a}
\shalf\dot{H}_{0} \leq -\nu H_{1} + H_{0}^{1/2}\|\bdf\|_{2}
\ee
in combination with Poincar\'{e}'s inequality $k_{1}^{2}H_{0} \leq H_{1}$, 
gives
\bel{Ler1b}
\dot{X} \leq -\nk X + \|\bdf\|_{2}
\ee
Let us also denote $X_{f}$ by ($k_{1} = 2\pi/L$)
\bel{en2}
X_{f} =  \frac{\|\bdf\|_{2}}{\nk} = \left(\frac{\nu}{L^{3/2}k_{1}^{2}}\right)\Gr
\leq c\,\nu L^{1/2}\Rey^{2}
\ee
which has the same dimensions as $X$. Integration of (\ref{Ler1b}) from $t_{0}$
to $t$ results in
\bel{en3}
X(t) \leq X_{0}\enk + X_{f}\left(1- \enk\right)
\ee
Because there is no specific knowledge of $t_{0}$ the upper bound on $H(t)$ is
taken over the full time-range $0\leq t\leq\infty$ which, from (\ref{en3}), is
\bel{en4}
H_{0}(t_{0}) \leq 
\left(\frac{\nu^{2}}{L^{3}k_{1}^{4}}\right)\Gr^{2}\leq c\,\nu^{2} L \Rey^{4}
\ee
This is properly valid after the time when transients have died out. 
The exponential decay in (\ref{en3}) is trivial compared to 
$\exp(\omega_{0}\Rey\,\Delta t)$ so we obtain 
\bel{en5}
\Idt e^{\omega_{0}\Rey\,\Delta t} X(t)\,dt \leq c\,\nu L^{1/2}
\calE(\Delta t )\Rey^{2}.
\ee
Now multiply (\ref{Ler1a}) by $e^{\omega_{0}\Rey\,\Delta t}$ and 
integrate by parts
\beq{en6}
\nu\Idt e^{\omega_{0}\Rey\,\Delta t}H_{1}(t)\,dt 
&\leq& \Idt e^{\omega_{0}\Rey\,\Delta t}\left(-\shalf\dot{H}_{0} + 
X\|\bdf\|_{2}\right)dt\nonumber\\
\leq \shalf H_{0}(t_{0}) - \shalf H_{0}(t)e^{\omega_{0}\Rey\,\Delta t}
&+& \shalf\omega_{0}\Rey\Idt H_{0} e^{\omega_{0}\Rey\,\Delta t}dt\nonumber\\
&+& c\,\nu^{3} L^{-1} \calE(\Delta t )\Rey^{4}
\eq
In the general case the negative term can be dropped, leaving
\beq{en8}
\nu\Idt e^{\omega_{0}\Rey\,\Delta t}H_{1}(t)\,dt &\leq& c_{1}\nu^{2}L\Rey^{4}\nonumber\\
&+& c_{2}\nu^{3}L^{-1}\calE(\Delta t )(\Rey^{5}+O(\Rey^{4}))
\eq
The predominant $\Rey^{5}$-term has a correction term of $O(\Rey^{4})$ from 
the fourth term in (\ref{en6}) and another of $O(\Rey^{2})$ from making 
up $H_{1}$ to $F_{1}$.  

Note that the first term on the right hand side of (\ref{F1est1}) in Lemma
\ref{lemma3} can be removed if the energy has the lower bound $H_{0}(t)\geq
H_{0}(t_{0})e^{-\omega_{0}\Rey\,\Delta t}$: see \S\ref{condreg} for a
discussion of this.
\end{proof}
\begin{lemma}\label{lemma4}
An estimate for the exponentially weighted time integral of $\kappa_{2,1}$ is
\beq{kappalem}
\Idt\kappa_{2,1}(t)e^{\omega_{0}\Rey\,\Delta t}\,dt 
&\leq& c\,\nu^{-2}\Idt e^{\omega_{0}\Rey\,\Delta t}F_{1}dt\nonumber\\
&+& c\,\calE(\Delta t )(L\Gr)^{-1}\Rey  
\eq
\end{lemma}
\begin{proof}
The time integral of $\kappa_{2,1}$ can be estimated from 
\beq{kappa5}
\Idt \kappa_{2,1}e^{\omega_{0}\Rey\,\Delta t}\,dt &= & 
\Idt e^{\omega_{0}\Rey\,\Delta t}(F_{2}/F_{1}^{2})^{1/2}F_{1}^{1/2}\,dt\nonumber\\
\leq  \frac{\,\nu^{2}}{2}\Idt e^{\omega_{0}\Rey\,\Delta t}(F_{2}/F_{1}^{2})\,dt
&+& \frac{1}{2\nu^{2}}\Idt e^{\omega_{0}\Rey\,\Delta t}F_{1}\,dt
\eq
The first integral on the far right hand side of (\ref{kappa5}) can be estimated
by using the inequality for $F_{1}$ from (\ref{ladderlemma2}) in Proposition 
\ref{Fnladder} in \S\ref{Fnsect}
\bel{kappa5a}
\shalf\dot{F}_{1} \leq -\frac{\nu}{4}F_{2} + c\,\nu^{-3}F_{1}^{3} + 
\omega_{0}\Rey F_{1}
\ee
Dividing (\ref{kappa5a}) by $F_{1}^2$, multiplying by $e^{\omega_{0}\Rey\,
\Delta t}$
and integrating 
\beq{kappa5b}
\frac{\nu}{4}\Idt e^{\omega_{0}\Rey\,\Delta t}
(F_{2}/F_{1}^{2})dt &\leq& c\,\nu^{-3}\Idt e^{\omega_{0}\Rey\,\Delta t}F_{1}dt\nonumber\\
& + & \shalf\left(F_{1}^{-1}(t)e^{\omega_{0}\Rey\,\Delta t} - F_{1}^{-1}(t_{0})\right)
\eq
The last term can be rewritten in terms of $\calE(\Delta t)$ which leaves a
$F_{1}^{-1}$ term. The upper bound on this is proportional to $\Gr^{-1}$, which
can be ignored as small.
\beq{kappa5c}
\Idt \kappa_{2,1}e^{\omega_{0}\Rey\,\Delta t}\,dt 
&\leq& c\,\nu^{-2}\Idt e^{\omega_{0}\Rey\,\Delta t}F_{1}dt\nonumber\\
&+& c\,L^{-1}\calE(\Delta t )\Gr^{-1}\Rey 
\eq
as in (\ref{kappalem}) above.
\end{proof}

\subsection{An estimate for $(\Delta t)_{b}$ when $n\geq 2$}\label{widthest}

The two estimates above for the weighted time integrals of $F_{1}$ and
$\kappa_{2,1}$ allow us to prove the main result of this section for 
$n\geq 2$. Define
\bel{epsilondef}
a_{n} = \frac{\lambda_{n+1}}{\mu}\left(\frac{2n-2}{2n-1}\right)-
\frac{4}{2n-1} - 5
\ee
Then $a_{n}\geq 1$ if $\mu$ is chosen such that
\bel{badwidth2a}
\mu \leq \lambda_{n+1}\left(\frac{n-1}{6n-1}\right)
\ee
whereas $0 < a_{n} < 1$ if $\mu$ is chosen to lie in the 
range\footnote{Note that for $n=2$ the lower bound on $\mu$ in 
(\ref{badwidth1a}) is greater than $\frac{1}{2(n+1)}$.}
\bel{badwidth1a}
\lambda_{n+1}\left(\frac{n-1}{6n-1}\right) < \mu <
\lambda_{n+1}\left(\frac{2n-2}{10n-1}\right)
\ee
\begin{theorem}\label{badwidththm}
For $n\geq 2$, if $a_{n} \geq 1$ the width of a bad interval is bounded by
\bel{badwidth2}
\tilde{c}_{n,1}\omega_{0}(\Delta t)_{b} \leq \Rey^{-a_{n}}
\ee 
whereas if $0 < a_{n} < 1$
\bel{badwidth1}
\tilde{c}_{n,2}\omega_{0}(\Delta t)_{b} \leq \Rey^{-1}\ln\Rey
\ee
\end{theorem}
\begin{remark}\label{n=1rem}
There appears to be no obvious parallel result for the finiteness of bad 
intervals in the case $n=1$, although Lemma \ref{n=1lemma} gives a lower 
bound $\mu > \shalf$ for the prevention of singularities forming on 
good intervals. 
\end{remark}
\begin{proof}
Let us return to Proposition \ref{Fnladder} in \S\ref{Fnsect}, 
inequality (\ref{ladderlemma4}), and recall that $\omega_{0} = \nu L^{-2}$
\bel{Yn1}
\shalf\dot{F}_{n} \leq - \shalf\nu F_{n+1} + c_{n}\left(\nu^{-1}
\|\bu\|_{\infty}^{2} + \omega_{0}\Rey\right)F_{n}.
\ee
This was manipulated in \S\ref{longav2} to produce (\ref{k11}) which is 
re-stated here as 
\bel{Yn3}
(n-\shalf)(\dot{Y}_{n} + \omega_{0}\Rey Y_{n}) \geq 
\shalf\nu\frac{F_{n+1}}{F_{n}^{\frac{2n}{2n-1}}} - c_{2}\nu^{-1}
\kappa_{2,1}^{\frac{1}{2n-1}}F_{1}^{\frac{2n-2}{2n-1}}
\ee
The first term on the right hand side of (\ref{Yn3}) can be estimated as
\bel{Yn5}
\frac{F_{n+1}}{F_{n}^{2n/(2n-1)}} \geq \kappa_{n+1}^{\frac{2n-2}{2n-1}}
F_{0}^{-\frac{1}{2n-1}}
\geq c\,\kappa_{n+1}^{\frac{2n-2}{2n-1}}(\nu^{2}L\Rey^{4})^{-\frac{1}{2n-1}}
\ee
having used the fact that $\kn \leq \kappa_{n+1}$. This result, together 
with a H\"{o}lder inequality, gives 
\beq{Yn6}
(n-\shalf)\frac{~d}{dt}\left[Y_{n}e^{\omega_{0}\Rey\,\Delta t}\right] 
&\geq & 
c\,\nu^{\frac{2n-3}{2n-1}}(L\Rey^{4})^{-\frac{1}{2n-1}}
e^{\omega_{0}\Rey\,\Delta t}\kappa_{n+1}^{\frac{2n-2}{2n-1}}\nonumber\\
&-& c_{2}^{\frac{2n-2}{2n-1}}e^{\omega_{0}\Rey\,\Delta t}
\left\{\nu^{\frac{2n-3}{2n-1}}\kappa_{2,1} + 
\nu^{-\frac{2n+1}{2n-1}}F_{1}\right\}
\eq
So far this has just been re-arrangement of (\ref{Yn1}). The lower bound 
on $\kn$ is now applied to the first term on the right hand side along with 
a time integration 
\beq{Yn7}
(n-\shalf)\left\{Y_{n}(t)e^{\omega_{0}\Rey\,\Delta t}\right\}
&\geq & c_{n+1}^{\frac{2n-2}{2n-1}}\nu^{\frac{2n-3}{2n-1}}L^{-1}\calE(\Delta t)
\Rey^{\frac{(2n-2)\lambda_{n+1}}{(2n-1)\mu}-\frac{4}{2n-1}}\nonumber\\
&-& c_{2}^{\frac{2n-2}{2n-1}}\nu^{\frac{2n-3}{2n-1}}\Idt e^{\omega_{0}\Rey\,\Delta t}\kappa_{2,1}\,dt\nonumber\\
&-& c_{2}^{\frac{2n-2}{2n-1}}\nu^{-\frac{2n+1}{2n-1}}\Idt e^{\omega_{0}\Rey\,\Delta t}F_{1}\,dt\nonumber\\
&+& (n-\shalf)Y_{n}(t_{0})
\eq
For the left hand side it is sufficient to show that this is bounded above by 
a very small number on a bad interval
\bel{Yn8}
Y_{n} = \kappa_{n+1}^{-\frac{2n}{2n-1}}F_{0}^{-\frac{1}{2n-1}}
\leq L\nu^{-\frac{2}{2n-1}}
\Rey^{-\frac{2n\lambda_{n+1}}{\mu(2n-1)}}\Gr^{-\frac{1}{2n-1}}
\ee
Using Lemmas \ref{lemma3} and \ref{lemma4} a comparison of the major terms in 
(\ref{Yn7}) shows that
\bel{Yn9}
\omega_{0}\calE(\Delta t)
\left\{c_{n+1}^{\frac{2n-2}{2n-1}}\Rey^{\frac{(2n-2)\lambda_{n+1}}{(2n-1)\mu}-
\frac{4}{2n-1}}-c_{3}\Rey^{5}\right\}
\leq c_{4}\Rey^{4}
\ee
For $n\geq 2$ the left hand side is always positive provided $\mu$ is chosen 
in the range
\bel{Yn10}
\frac{\lambda_{n+1}}{\mu}\left(\frac{2n-2}{2n-1}\right)-\frac{4}{2n-1} > 5
\ee
or in the range
\bel{Yn11}
\mu < \lambda_{n+1}\left(\frac{2n-2}{10n-1}\right)
\ee
To solve (\ref{Yn9}) use the definition of $a_{n}$ in (\ref{epsilondef}) 
giving
\bel{Yn13}
\tilde{c}_{n}\,\omega_{0}\calE(\Delta t)\leq\Rey^{-a_{n}}
\ee
The solution of this depends on whether $a_{n}$ lies in the range 
$a_{n}\geq 1$ or $0< a_{n} <1$. The estimates in (\ref{Einequal1}) 
and (\ref{Einequal2}) are appropriate.
\end{proof}


\subsection{Intersection of bad intervals: the relation to Scheffer's 
singular set}\label{badintersect}

Figure \ref{knfig} of \S2.4 is a representation of good and bad intervals 
for some $n\geq 2$. Since it must be assumed that the position of the 
intervals changes with $n$, the intersection of all the bad intervals 
for $n\geq 2$ is pertinent: only if this intersection is non-empty 
will singularities be possible. For each $n \geq 2$, let us 
designate a bad interval as the set $\mathcal{B}_{n}$ on the time-axis 
on which 
\bel{in4}
c_{n}\frac{\kappa_{n+1}}{\kn} <(L\kn)^{\mu}\Rey^{-\lbd_{n}}
\ee
Moreover, because $L\kn\leq L\kappa_{n+1}$, on this set there is a lower 
bound
\bel{lb1}
(L\kn)^{\mu} > c_{n}\Rey^{\lambda_{n}}.
\ee
Now consider the set $\mathcal{B}_{n+1}$ on which 
\bel{in5}
c_{n+1}\frac{\kappa_{n+2}}{\kappa_{n+1}} < 
(L\kappa_{n+1})^{\mu}\Rey^{-\lbd_{n+1}}
\hspace{.5cm}\Rightarrow\hspace{.5cm}
(L\kappa_{n+1})^{\mu} > c_{n+1}\Rey^{\lambda_{n+1}}
\ee 
Then on the intersection 
$\mathcal{I}_{n+1} = \mathcal{B}_{n}\cap\mathcal{B}_{n+1},$
we have 
\bel{in7}
c_{n}^{1+\mu}c_{n+1} (L\kappa_{n+2}) < (L\kn)^{\mpo^{2}}
\Rey^{-\lbd_{n+1}- \mpo\lbd_{n}}
\ee
Using $(L\kappa_{n+2})^{\mu} \geq (L\kappa_{n+1})^{\mu} > 
c_{n+1}\Rey^{\lambda_{n+1}}$ on $\mathcal{I}_{n+1}$, a new lower bound is
\bel{in8}
(L\kn)^{1+\mu} > \left(c_{n+1}\Rey^{\lbd_{n+1}}\right)^{1/\mu}
\left(c_{n}\Rey^{\lbd_{n}}\right)
\ee
Now consider the intersection 
$\mathcal{I}_{n+2} = \mathcal{B}_{n}\cap\mathcal{B}_{n+1}\cap\mathcal{B}_{n+2}$.
On this set there is a larger lower bound
\bel{in10}
(L\kn)^{(1+\mu)^{2}} > \left(c_{n+2}\Rey^{\lbd_{n+2}}\right)^{1/\mu}
\left(c_{n+1}\Rey^{\lbd_{n+1}}\right)\left(c_{n}\Rey^{\lbd_{n}}\right)^{1+\mu}
\ee
We wish to find a lower bound on $L\kn$ on the set of $p$ intersections
\bel{pinter}
\mathcal{I}_{n+p} = \mathcal{B}_{n}
\cap\mathcal{B}_{n+1}\cap\ldots\cap\mathcal{B}_{n+p}
\ee
By inspection, the general formula for the lower bound of $L\kn$ 
on $\mathcal{I}_{n+p}$ is
\begin{eqnarray}\label{g1}
(L\kn)^{(1+\mu)^{p+1}} &> &
\left(c_{n+p+1}\Rey^{\lbd_{n+p+1}}\right)^{1/\mu}
\left(c_{n+p}\Rey^{\lbd_{n+p}}\right)
\ldots\left(c_{n}\Rey^{\lbd_{n}}\right)^{(1+\mu)^p}\nonumber\\
&=& \left(c_{n+p+1}\Rey^{\lbd_{n+p+1}}\right)^{1/\mu}
\left(\Pi_{i=0}^{p}c_{n+i}^{\xi_{p,i}}\right)\Rey^{L_{n,p}}
\end{eqnarray}
where
\bel{biglambdadef}
\xi_{p,i} = \mpo^{p-i}
\hspace{1.5cm}
L_{n,p} = \sum_{i=0}^{p}\lbd_{n+i}\,\xi_{p,i}
\ee
We are particularly interested in the limit $p\to\infty$ so we write 
\bel{g4}
(L\kn)^{\mu} > \Rey^{\Lambda_{n}^{(\infty)}}
\ee
where $\Lambda_{n}^{(p)}$ is defined as 
\bel{g4a}
c_{n}^{(\infty)}\Rey^{\Lambda_{n}^{(p)}} = 
\left\{\left(c_{n+p+1}\Rey^{\lbd_{n+p+1}+\mu\,L_{n,p}}\right)
\left(\Pi_{i=0}^{p}c_{n+i}^{\xi_{p,i}}\right)^{\mu}
\right\}^{\frac{1}{(1+\mu)^{p+1}}}
\ee
and then, because $\lbd_{n+p} > \lbd_{n}$ and $c_{n+p} > c_{n}$, it follows 
that
\bel{g3}
\xi_{p}\lambda_{n} < L_{n,p} < \xi_{p}\lambda_{n+p}
\hspace{1.5cm}
c_{n}^{\,\xi_{p}}< \Pi_{i=0}^{p}\,c_{n+i}^{\xi_{p,i}}
< c_{n+p}^{\,\xi_{p}}
\ee
where $\xi_{p}$ is defined by the sum 
\bel{g2}
\xi_{p} = \sum_{i=0}^{p}\xi_{p,i} = \mu^{-1}\left\{(1+\mu)^{p+1}-1\right\}
\ee
\par\vspace{-2cm}
\newcommand{\Z}{\rule[0.0cm]{0.36cm}{0.1cm}} 
$$
\begin{minipage}[ht]{9cm}\label{figstrip} 
\setlength{\unitlength}{.75cm}
\begin{picture}(11,11)(0,0)
\thicklines
\put(0,0){\vector(0,1){8}}
\put(0,0){\vector(1,0){10}}
\put(0,8.2){\makebox(0,0)[b]{$\kappa_{n}(t)$}}
\put(10.2,0){\makebox(0,0)[b]{$t$}}
\thinlines
\put(3,0){\line(0,1){8}}
\put(3.5,0){\line(0,1){8}}
\put(7,0){\line(0,1){8}}
\put(7.5,0){\line(0,1){8}}
\multiput(0,2)(.1,0){100}{.}
\multiput(0,6)(.1,0){30}{.}
\put(3,6){\Z}
\multiput(3.5,6)(.1,0){65}{.}
\put(-.75,2){\makebox(0,0)[b]{\scriptsize$Re^{\lambda_{n}}$}}
\put(5.3,2.1){\makebox(0,0)[b]{\scriptsize Long-time average}}
\put(5.2,6.5){\makebox(0,0)[b]
{\scriptsize$(L\kappa_{n})^{\mu}>Re^{\lambda_{n}}$}}
\put(3.9,6.4){\vector(-1,0){.65}}
\put(6.2,6.4){\vector(1,0){1.1}}
\put(-.75,6){\makebox(0,0)[b]{\scriptsize$Re^{\lambda_{n}/\mu}$}}
\put(3.25,-.70){\makebox(0,0)[b]{\scriptsize$(\Delta t)_{b}$}}
\put(5.25,-.70){\makebox(0,0)[b]{\scriptsize$(\Delta t)_{g}$}}
\put(5.75,-.25){\vector(1,0){1.25}}
\put(4.75,-.25){\vector(-1,0){1.2}}
\qbezier[250](0,1)(3,0)(3.01,8)
\qbezier(3.5,7.75)(5,-7)(7,7.75)
\qbezier[300](7.5,8)(7.5,0)(10,1)
\end{picture}
\end{minipage}
$$
\par\vspace{.5cm}\noindent
{\small \textbf{Figure 3:} Similar to figure \ref{knfig}, representation
of good/bad intervals for $\kappa_{n}$ with a black strip representing 
the bad interval used in the intersection table 
\ref{intersect1}.} 
\par\medskip\noindent
The potentially singular set, $\mathcal{S}^{(\infty)}$, is given by 
\bel{Iinfty}
\mathcal{S}^{(\infty)} = \mathcal{B}_{1}\cap\mathcal{B}_{2}\cap\ldots
\cap\mathcal{B}_{n}\cap\ldots
\ee
must necessarily include $\mathcal{B}_{1}$, the singular set of $\kappa_{1}$
(and therefore $F_{1}$). The range of values of $\mu$ expressed in
(\ref{badwidth1a}) and Theorem \ref{badwidththm}) are valid for $n\geq 2$.
These narrow to $0 < \mu < \thfth$ in the limit $n\to\infty$. As already
pointed out in Lemma \ref{n=1lemma} of \S\ref{good}, a corresponding separate
calculation for $F_{1}$ shows that $\mu$ lies in the range $\shalf < \mu < 1$
for $n=1$. When the allowed ranges of $\mu$ are taken into account for good
and bad intervals we conclude that 
\begin{theorem}\label{muthm} For all bad
intervals to be finite for $n\geq 2$ and for no singularities to form in good
intervals for $n\geq 1$ the allowed range of $\mu$ is
\bel{murange}
\shalf < \mu < \thfth\,.
\ee
\end{theorem}
$\mathcal{S}^{(\infty)}$ is related to Scheffer's potentially singular set: 
his set is technically the union of all sets $\mathcal{S}^{(\infty)}$ 
associated with \textit{every} bad interval. Scheffer showed that 
this set has zero $\shalf$-dimensional Hausdorff measure \cite{Scheffer76a}, 
which means that it must consist of, at most, points. Whether the $\kn$ 
actually become singular on this set is still an open question. 
From (\ref{g4a})-(\ref{g2}) we have
\bel{g5}
c_{n}\Rey^{\lbd_{n}}< c_{n}^{(\infty)}\Rey^{\Lambda_{n}^{(\infty)}} 
<\lim_{p\to\infty}c_{n+p+1}\Rey^{\lbd_{n+p+1}}
\ee
Divergence in this limit would guarantee singular behaviour if the set
$\mathcal{S}^{(\infty)}$ is non-empty but there is no evidence that the
product in (\ref{g4a}) diverges in the limit even though the upper bound in
(\ref{g5}) is infinite. From (\ref{g5}) it is clear that
$\Lambda_{n}^{(\infty)}>\lambda_{n}$ so all the estimates of the previous
sections dependent upon $\lambda_{n}$ should be replaced by
$\Lambda_{n}^{(\infty)}$. This paper, however, furnishes no evidence on 
the distribution of the intervals; Table 6 is simply a pictorial 
representation of some randomly chosen bad intervals associated with 
$\kappa_{n}\to\kappa_{n+8}$ to illustrate how the final intersection 
may form.
\arrayrulewidth=0.3mm 
\doublerulesep=0.3mm 
\begin{table}[t] 
\begin{center} 
\begin{tabular}{||c||@{\hspace{0.mm}}c@{\hspace{0.mm}}
                    @{\hspace{0.mm}}c@{\hspace{0.mm}}
                    @{\hspace{0.mm}}c@{\hspace{0.mm}}
                   @{\hspace{0.mm}}c@{\hspace{0.mm}} 
                    @{\hspace{0.mm}}c@{\hspace{0.mm}} 
                    @{\hspace{0.mm}}c@{\hspace{0.mm}} 
                   @{\hspace{0.mm}}c@{\hspace{0.mm}} 
                    @{\hspace{0.mm}}c@{\hspace{0.mm}} 
                    @{\hspace{0.mm}}c@{\hspace{0.mm}}|| 
                   @{\hspace{0.mm}}c@{\hspace{0.mm}}
                    @{\hspace{0.mm}}c@{\hspace{0.mm}}
                    @{\hspace{0.mm}}c@{\hspace{0.mm}}|} 
\hline 
$\mathcal{I}_{n+8}$&\Q\Q\Q&\Q\Q\Q&\Y\Q\Q&\Q\Q\Q      
                       &\Q\Y\Q&\Q\Q\Q&\Q\Y\Q&\Y\Q\Q&\Q\Q\Y \\
\hline 
\hline 
$\kappa_{n+8}$                &\X\X\X&\Q\Q\X&\X\Q\Q&\X\X\X 
                      &\X\X\X&\Q\Q\Q&\X\X\Q&\X\X\X&\X\X\X \\
$\kappa_{n+7}$                &\X\Q\Q&\X\X\X&\X\Q\X&\X\X\X 
                      &\X\X\X&\Q\Q\X&\X\X\X&\X\Q\Q&\Q\X\X \\
$\kappa_{n+6}$                &\X\X\Q&\Q\X\X&\X\X\X&\X\X\X 
                      &\X\X\X&\Q\X\X&\X\X\X&\X\X\Q&\Q\Q\X \\
$\kappa_{n+5}$                &\X\X\Q&\X\X\X&\X\X\X&\Q\Q\X 
                      &\X\X\X&\X\X\Q&\Q\X\X&\X\X\X&\X\X\X \\
$\kappa_{n+4}$                &\Q\X\X&\X\Q\Q&\X\X\X&\X\X\X 
                      &\X\X\X&\Q\X\X&\X\X\X&\X\Q\X&\X\X\X \\
$\kappa_{n+3}$                &\X\X\Q&\Q\Q\X&\X\X\X&\X\Q\Q 
                      &\X\X\X&\X\X\X&\X\X\X&\X\X\X&\Q\X\X \\
$\kappa_{n+2}$                &\X\X\X&\X\X\X&\X\X\X&\X\Q\Q 
                      &\Q\X\X&\X\X\X&\X\X\X&\X\X\X&\Q\Q\X \\
$\kappa_{n+1}$               &\X\X\X&\X\X\X&\X\X\X&\Q\Q\X 
                      &\X\X\Q&\Q\X\Q&\X\X\X&\X\X\X&\X\X\X \\
$\kappa_{n}$                 &\X\X\X&\X\X\X&\X\X\X&\X\X\X 
                      &\X\X\X&\X\X\X&\X\X\X&\X\X\X&\X\X\X \\
\hline 
\end{tabular} 
\end{center} 
\caption{The lowest continuous horizontal black strip is the bad interval 
of $\kappa_{n}$ shown as the black strip in Figure 3. The strips in next 
8 levels are an illustration of how some randomly chosen bad intervals 
($\kappa_{n}\to\kappa_{n+8}$) could appear. The thicker strips at the 
highest level are the intersection of the 9 strips below.}\label{intersect1} 
\end{table} 

\subsection{Dangerous sub-intervals}\label{subint}

In addition to the intersection idea of the last section, we consider the 
special set of sub-intervals within each bad interval on which $\dot{F}_{n}
\geq 0$. Consider the $j$th sub-interval within the $i$th bad interval: this
is designated as \textit{dangerous} sub-interval of width $\tplus$ with an
initial value of designated as $\tij$. It is on these where singularities are
possible: they are not possible where any one of the $F_{n}$ is decreasing.
(\ref{badwidth2a}) and (\ref{badwidth1a}) show that the smallest lower bound 
on $\lambda_{n}/\mu$ is $\lambda_{n}/\mu >5$. Because $\lambda_{n+1} > 
\lambda_{n}$, and replacing $\lambda_{n+1}$ by $\Lambda_{n+1}^{(\infty)}$, 
define  
\bel{bndef}
b_{n} = \frac{\Lambda^{(\infty)}_{n+1}}{\mu} - 4 > 1
\ee
\begin{theorem}\label{subthm}
Dangerous sub-intervals are bounded in width by ($n\geq 2$)
\bel{sub0a}
\omega_{0}\tplus \leq c_{n} \Rey^{-b_{n}}\hspace{1.75cm}b_{n}> 1
\ee
and on these sub-intervals
\bel{sub0c}
L^{-3}F_{1}\geq \omega_{0}^{2}\Rey^{b_{n}+4}
\hspace{.75cm}\|\bu\|_{\infty}\geq L\omega_{0}\Rey^{b_{n}+4}
\ee
\bel{sub0d}
\|\nabla\bu\|_{\infty}\geq \omega_{0}\Rey^{b_{n}+4} 
\ee
\end{theorem}
\begin{remark}
Note that $b_{n} > a_{n}$ so the upper bounds of these sub-intervals 
are smaller than those in Theorem \ref{badwidththm}. 
\end{remark}
\begin{proof}
Consider (\ref{kladder0})
\bel{sub1}
\shalf\dot{F}_{n} \leq -\shalf\nu F_{n+1} + 
\left(c_{n}\nu^{-1}\|\bu\|^{2}_{\infty} + \omega_{0}^{2}\Rey\right)F_{n}.
\ee
Now use the Sobolev inequality $\|\bu\|_{\infty}^{2} \leq c\,
\kappa_{n,1}F_{1}$, and divide (\ref{sub1}) by $F_{n}$. Then on 
these sub-intervals
\bel{sub2}
\kappa_{n+1,n}^{2} \leq c_{n}\nu^{-2}\kappa_{n,1}F_{1} +  c\,L^{-2}\Rey
\ee
Now we know that $\kappa_{n+1,n} \geq \kappa_{n+1}$ so
\bel{sub3}
\kappa_{n+1}^{2} \leq \left(c_{n}\nu^{-2}F_{1}\right)^{2} +  c\,L^{-2}\Rey
\ee
Now the lower bound $L\kappa_{n+1}\geq \Rey^{\Lambda^{(\infty)}_{n+1}/\mu}$ 
is invoked giving 
\bel{sub4}
\left(c_{n}^{(\infty)}\Rey^{\Lambda^{(\infty)}_{n+1}/\mu}\right)^{2} - 
c\,\Rey \leq \left(c_{n}L\nu^{-2}F_{1}\right)^{2}
\ee
Because $\lambda_{n}/\mu >1$ for $n\geq 2$, the $\Rey$ term is small in
comparison, leaving 
\bel{sub5}
c_{n}^{(\infty)}\Rey^{\Lambda^{(\infty)}_{n+1}/\mu} \leq L\nu^{-2}F_{1}
\ee
Multiplying by the exponential term, integrating over $\Delta t^{+}$,  
and then using the exponentially time-weighted integral of $F_{1}$ in 
Lemma \ref{lemma3} gives
\bel{sub6}
c_{n}\,\omega_{0}\calE(\Delta t^{+})
\left[\Rey^{\Lambda^{(\infty)}_{n+1}/\mu}-\Rey^{5}\right]\leq \Rey^{4}
\ee
Now we know that $\Lambda^{(\infty)}_{n+1}/\mu >5$ and 
$b_{n}>1$ so the result in (\ref{sub0a}) follows. The definition of 
$a_{n}$ in (\ref{epsilondef}) guarantees that 
\bel{sub7}
b_{n} = \frac{\Lambda^{(\infty)}_{n+1}}{\mu} - 4 > a_{n} 
\ee
which is the correct way round. (\ref{sub0c}) and (\ref{sub0d}) 
follow from (\ref{sub5}).
\end{proof}

\section{A conditional regularity result}\label{condreg}

The reader who has followed the proof of Theorem \ref{badwidththm} will have
noticed that attempts to prove Navier-Stokes regularity fail in the bad
intervals. There, use was made of the variables $Y_{n}(t) = F_{n}^{-1/(2n-1)}$
defined in (\ref{Yndef}). To prevent the formation of singularities, it would
be necessary to show that $Y_{n}$ can never touch zero in a finite time.
Within the dangerous sub-intervals of \S\ref{subint} this can be achieved
provided the energy is bounded below in a certain manner. Specifically the
result is
\begin{theorem}\label{condregthm}
The Navier-Stokes equations are regular if, in dangerous sub-intervals $\tplus$, 
there is a lower bound on the energy 
\bel{cond1}
H_{0}(t)\geq H_{0}(\tij)e^{-\omega_{0}\Rey\Delta t}
\ee
\end{theorem}
\begin{remark}
Over the very short time interval $\Delta t$ the exponent on the right 
hand side of (\ref{cond1}) is very small, so the right hand side 
is almost $H_{0}(\tij)$.
\end{remark}
\begin{proof}
They key point preventing progress regarding bounding
$Y_{n}(t)$ away from zero is the set of extra terms in Lemmas \ref{lemma3} and
\ref{lemma4} that are not coefficients of $\calE (\Delta t)$. This creates
negative terms on the right hand side of (\ref{Yn7}) that cannot be controlled. 
To circumvent this problem it is necessary to remove two hurdles. The first is
the last pair of terms within inequality (\ref{kappa5b})
\beq{cond2}
\frac{\nu}{4}\Idt e^{\omega_{0}\Rey\,\Delta t}
(F_{2}/F_{1}^{2})dt &\leq& c\,\nu^{-3}
\Idt e^{\omega_{0}\Rey\,\Delta t}F_{1}dt\nonumber\\
& + & \shalf\left[F_{1}^{-1}(t)e^{\omega_{0}\Rey\,\Delta t} 
- F_{1}^{-1}(t_{0})\right]
\eq
On dangerous sub-intervals $\tij$ where $F_{1}$ is increasing, the last term
in (\ref{cond2}) can be re-written as
\bel{cond3}
\shalf\left(F_{1}^{-1}(t)e^{\omega_{0}\Rey\,\Delta t} - F_{1}^{-1}(\tij)\right)
\leq \shalf \omega_{0}\Rey\,F_{1}^{-1}(\tij)\calE\tplus
\ee
This term is now classed as one of the $\calE\tplus$ terms and is merely a
term of lower order than the dominant $\Rey^{5}$ term.

Secondly, in (\ref{en6}) if it assumed that on these sub-intervals that 
(\ref{cond1}) is true then the extra terms in (\ref{F1est1}) can be removed, 
leaving
\bel{cond4}
\Idt e^{\omega_{0}\Rey\,\Delta t}F_{1}(t)\,dt \leq  
c_{2}\nu^{2}L^{-1}\calE\tplus \left[\Rey^{5}+O(\Rey^{4})\right]
\ee
which again is proportional to $\calE\tplus$. Thus (\ref{Yn7}) becomes 
\beq{cond5}
(n-\shalf)\left\{Y_{n}(t)e^{\omega_{0}\Rey\,\tplus}\right\}
&\geq & (n-\shalf)Y_{n}(\tij)\\
&+& c_{n}\nu^{\frac{2n-3}{2n-1}}L^{-1}\calE\tplus
\left\{\Rey^{\frac{(2n-2)\lambda_{n+1}-4\mu}{(2n-1)\mu}}-\Rey^{5}\right\}
\nonumber
\eq
with no negative terms on the right hand side. Given that $\mu$ is chosen in 
the restricted ranges in (\ref{badwidth2a}) and (\ref{badwidth1a}) within 
Theorem \ref{badwidththm}, and that $\calE >0$ for $t>0$, then $Y_{n}(t)$ 
can never be zero.
\end{proof}

\section{Discussion}\label{disc}

To summarize the arguments of this paper, it has been shown that very strong 
fluctuations in the $\kn(t)$ can occur in time, lower bounds on which are 
much higher than the long-time average (\ref{I6})
\bel{I6a}
\left<L\kn\right> \leq c_{n}\Rey^{\lambda_{n}}\hspace{2cm}
\lambda_{n}= 3-\frac{5}{2n} +\frac{\delta}{n}
\ee
This is based on the raising of the lower bound on the ratio $\kappa_{n+1}/\kn$ 
away for unity, a result which is expressed in Theorem \ref{intervalthm} 
of \S\ref{inter}
\bel{newinequal}
c_{n}\frac{\kappa_{n+1}}{\kappa_{n}}\geq (L\kn)^{\mu}\Rey^{-\lambda_{n}}
\ee 
and which is effective only on the good parts of the time axis. On those 
parts of the time-axis where the reverse of (\ref{newinequal}) is true, 
no upper bound has been found on the $\kappa_{n}$ but very large lower 
bounds exist on these of the form
\bel{newlb}
L\kappa_{n} > c_{n}\Rey^{\lambda_{n}/\mu}
\ee
The above results are valid for $n\geq 2$. By including intervals at $n=1$ 
the intersection set of bad intervals $\mathcal{S}^{(\infty)}$ is related 
to Scheffer's singular set of potential singularities, in which case the 
right hand side of (\ref{newlb}) can be raised again by replacing 
$\lambda_{n}$ by $\Lambda^{(\infty)}_{n}$. The constant\footnote{In fact 
(\ref{newinequal}) breaks the dimensional scaling of the standard Sobolev 
and Gagliardo-Nirenberg inequalities although how the introduction of the 
exponent $\mu$ affects this in a precise manner is not yet clear.} $\mu$ 
is then constrained to the range $\shalf <\mu<\thfth$. 

A picture emerges of Navier-Stokes solutions that are regular on `most' 
of the time-axis which is punctured by short, active intervals. While no 
upper bound on $\kappa_{n}$ has yet been found within these intervals, 
to become singular $\kappa_{n}$ would have to find its way through a 
non-empty intersection in a similar manner as the illustration in Figure 
(\ref{figstrip}). Notwithstanding the credence that must be given to the 
long-standing and widely held belief that the Navier-Stokes equations 
must be regular for arbitrarily long times, an alternative emerges that 
no upper bound exists on the $\kappa_{n}$ but the potentially singular 
set $\mathcal{S}^{(\infty)}$ in \S\ref{badintersect} is either empty or 
only allows extremely rare singular events. 
 
The results in this paper are consistent with ideas that have existed for 
many decades concerning intermittent flows \cite{BT49,Kuo71,SM88,MS91,Homoturb} 
but it has to be acknowledged that our results are lacking in four areas:
\begin{enumerate}
\item The bounds are almost certainly lacking in sharpness; the two main 
places where the size of the bounds appears is in $\lambda_{n}$ and in 
the $\Rey^{4}$ global estimate for $H_{0}$ in (\ref{en4}), which is 
certainly not sharp. Given these results, the state of the analysis 
is such that it may be premature to suggest specific numerical tests.

\item In addition to a lack of control over the $\kn$ within the bad 
intervals, their distribution and sensitivity to the value of $n$ is 
an important but unanswered question. 

\item The nature of solutions in the good  intervals has yet to be properly 
established. Because solutions are bounded pointwise in time and are also 
constrained by the long-time average (\ref{I6a}) it is to be expected that 
they should show a strong degree of quiescence, particularly within the 
central parts of these intervals. This has yet to be demonstrated. 

\item How solutions at the junctions of the good and bad intervals 
connect to each other is not clear as we have only weak solutions 
within the bad.
\end{enumerate}
Results of this type derived in the manner of \S\ref{inter} are not confined 
to the three-dimensional Navier-Stokes equations but could be applied to 
simpler problems. All that is needed are long-time average bounds for weak
solutions constructed in such a way that the equivalents of $\kn$ are bounded
below. An example that springs to mind is the case of the two-dimensional
Navier-Stokes equations. Not only are these regular but tight estimates exist
both for the attractor dimension \cite{CFT} and the number of determining
modes and nodes \cite{JT}. Other examples might be the alpha and Leray models
of turbulence \cite{FHT1,FHT2,CHOT} or the complex Ginzburg-Landau equation
\cite{CGL1,CGL2}.
 
The physics community has for many years used the scaling arguments based on
Kolmogorov's original work. Frisch's book gives a detailed factual and
historical account of these arguments \cite{Frisch}. Using an inertial range
argument, it has been argued in \cite{Doering02} that the scaling in the
rigorous upper bounds on $\left<\kn\right>$ from (\ref{I6}) (repeated in
(\ref{I6a}) above) may be interpreted in terms of the structure of the Fourier
spectrum $E_{s}(k)$ if a scaling of the form $E_{s}\sim k^{-q}$ is assumed in
the inertial range, up to the cut-off wave-number $L k_{c}\sim \Rey^{q_{c}}$.
Disregarding the correction from $\delta$, the \textit{a priori} bounds in
(\ref{I6a}) are consistent with $q = 8/3$ and $q_{c} = 3$. Such a $k^{-8/3}$
spectrum has arisen in at least two previous studies. Sulem \& Frisch
\cite{SF} have shown that a $k^{-8/3}$ spectrum is the borderline steepness
capable of sustaining an energy cascade in the Navier-Stokes equations when
the total energy is finite. Mandelbrot \cite{Man1}, and later Frisch, Sulem \&
Nelkin in their toy $\beta$-model \cite{FSN,Frisch}, came upon this same
scaling exponent as an extreme limit of intermittency in the energy cascade.
They found that if the energy dissipation is assumed to be concentrated on a
fractal set (in space) of dimension $D = 8-3q$, then the energy spectrum
scaling is of the form $E_{s}\sim k^{-q}$. Within this picture, the exponent
$q=8/3$ thus corresponds to dissipation concentrated at zero-dimensional
points in space. Interestingly, the conventional Kolmogorov $k^{-5/3}$
spectrum for homogeneous isotropic turbulence is associated with $D=3$; that
is, a complete lack of intermittency with dissipation spread uniformly in
space is consistent with $q = 5/3$. Departures from Kolmogorov scaling,
otherwise known as anomalous scaling, can be associated with intermittency in
the inertial-range. These arguments have been applied to and tested on various
models such as the $\beta$-model \cite{FSN}, and the bifractal and
multi-fractal models \cite{Frisch}. While these suffer by comparison in not
having the same degree of complexity as the Navier-Stokes equations -- as
Frisch \cite{Frisch} and Sreenivasan \cite{Sreeni85} have both pointed out --
these models have the merit of simplicity while capturing the main essence of
the phenomena in question. More recent work on anomalous scaling has centred
on the role of the $SO(3)$ symmetry group in the expansion of the correlation
functions \cite{ALP,ABMP}. \begin{acknowledgement} Thanks are due to Daryl
Hurst and Christos Vassilicos of Imperial College London Aeronautics
Department for the data in Figure 1 and to Robert Kerr of Warwick University,
and to Jean-Luc Thiffeault and Darryl Holm of Imperial College London
Mathematics Department for discussions. This work was begun at the 1999
Turbulence Program at the Isaac Newton Institute Cambridge and continued while
JDG was a Visiting Professor at RIMS, Kyoto University, in the summer of 2000.
\end{acknowledgement} 

\address{Department of Mathematics, Imperial College London, London SW7 2AZ UK 
\and Department of Mathematics \& Michigan Center for Theoretical Physics,
University of Michigan, Ann Arbor, Michigan,\\MI 48109-1109, USA}

\end{document}